\documentclass{article}
\usepackage{amsfonts}
\usepackage{graphicx}
\def\bs{\bigskip}  
\def\semi{\hbox{ $\times $ \kern-.972em \raise.12719em\hbox{ $_{^|}$}  }}
\def\reals{\mathbb R}
\def\ints{\mathbb Z}

\newtheorem{theorem}{Theorem}
\newtheorem{lemma}{Lemma}[section]
\newtheorem{corollary}{Corollary}
\def\pf {{\bf Proof:} \ }
\def\endpf{$\|$ \bigskip}
\def\be{\begin{enumerate}}
\def\ee{\end{enumerate}}
\def\bi{\begin{itemize}}
\def\ei{\end{itemize}}

\begin{document}
\title{OBSTRUCTIONS TO TRIVIALIZING A KNOT}
\author{Joan S. Birman and John Atwell Moody}
\date{submitted June 11, 2002; revised October 16, 2003}
\maketitle

\begin{abstract} The recent proof by Bigelow and Krammer that the braid
groups are linear opens the possibility of applications to
the study of knots and links. It was proved by the first author and Menasco
that
any closed braid representative of the unknot can be systematically simplified
to a round planar circle by a finite sequence of exchange moves and reducing
moves.  In this paper we establish connections between the faithfulness of the
Krammer-Lawrence representation and
the problem of recognizing when the conjugacy class of a closed braid admits an
exchange move or a reducing move.
\end{abstract}

\flushleft

\section{Introduction}
\label{section:introduction}
The goal of this article
is to study some problems in algorithmic knot theory with the use of new tools.
Our
work will combine two threads of thought:
\bi
\item  The first author's work with Menasco on the study of knots via
closed braids.
\item  The work of the second author, who proved that the faithfulness of
certain matrix representations of the braid group rested on
whether the matrices in question effectively detected intersections between
certain
arcs on the punctured plane and their images under braid homeomorphisms.
\ei

Let ${\cal K}$ be an oriented knot type in oriented 3-space $S^3$, or
alternatively in $\reals^3$.  A representative $\tilde{K} \in \cal{K}$ is said
to be
a {\em closed braid} if there is an unknotted curve ${\bf A}$ (the {\em braid
axis}), which we shall think of as the z-axis in $\reals^3$, and a
choice of fibration of $\reals^3 - {\bf A}$  by half-planes
$\{P\times\{t\}, t\in [0,1]\}$,  such that whenever $\tilde{K}$ meets a fiber
$P\times \{t\}$ the intersection is transverse. This implies that $\tilde{K}
\subset (\reals^3-{\bf A})$. The {\em braid
index}
$n = n(\tilde{K})$ is the number of points in $\tilde{K}\cap
(P\times\{t\})$. This number is independent of the choice of $t\in[0,1]$,
because
$\tilde{K}$ is transverse to every fiber $P\times\{t\}$.  Two closed braids
$\tilde{K},
\tilde{K'}\subset (\reals^3-{\bf A})$ are equivalent if there is is an isotopy
$\varphi:(\reals^3-{\bf A})\times[0,1]\to (\reals^3-{\bf A})$ with
$\varphi(\tilde{K},0) =\tilde{K}$ and $\varphi(\tilde{K},1)=\tilde{K'}$ and
each
$\varphi(\tilde{K},s)$ a closed braid.
Cutting $\reals^3 - {\bf A}$  along any half-plane $P\times\{t_0\}$  we
obtain a {\it braid} $K$.  Two open braids obtained by cutting along the same
half-plane
$P\times\{t_0\}$ are equivalent if there is an isotopy $\varphi$ as above which
is the
identity on the cutting plane. \bs

Open braids are in 1-1 correspondence with elements in Artin's braid group
${\bf
B}_n$ and closed braids are in 1-1 correspondence with conjugacy classes in
${\bf B_n}$. The open braids obtained by cutting along 
distinct planes are conjugate in
${\bf B}_n$. \bs

\begin{figure}[htpb!]
\centerline{\includegraphics[scale=.8, bb=135 232 487 680]{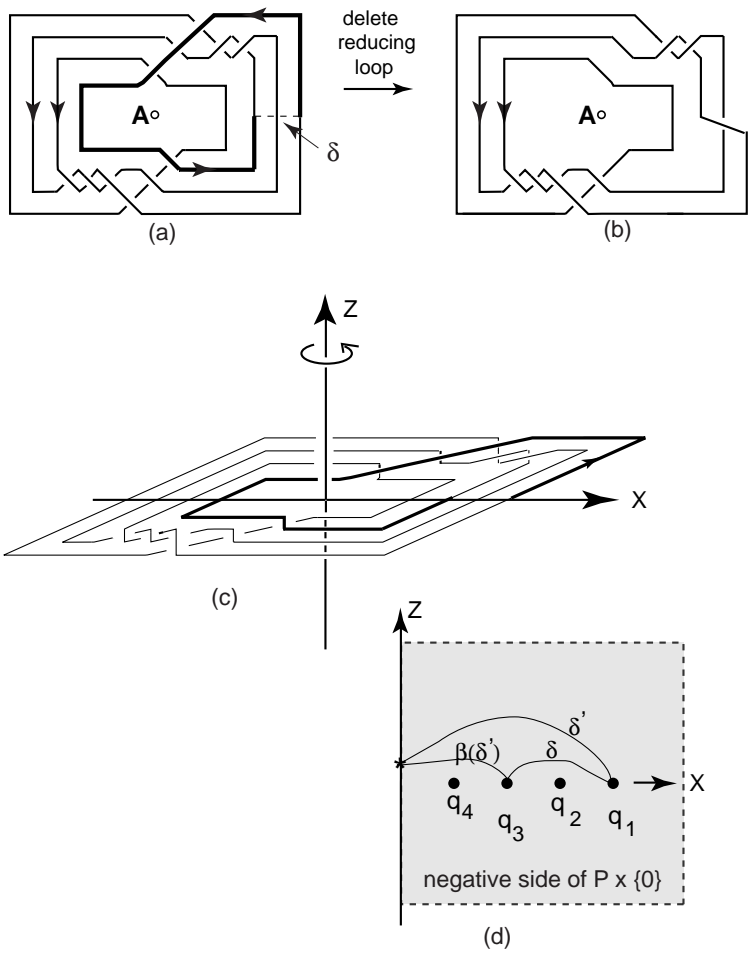}}
\caption{The closed 4-braid $\sigma_2^{-2}
\sigma_1^{-1}\sigma_2^{-1}\sigma_3^{-1}\sigma_2^3\sigma_1\sigma_2 \sigma_3$
 has a reducing loop. } 
\label{figure:reducible1}
\end{figure}

We now define what we mean when we say that the conjugacy class of a
closed braid $\tilde{K}$ `admits a reducing move'. An example is given in Figure \ref{figure:reducible1}. 
The plane $P\times
\{0\}$ meets $\tilde{K}$ in $n$ points $q_1,\dots,q_n$. These n points divide
$\tilde{K}$ into $n$ arcs $\alpha_1,\dots,\alpha_n$, where each $\alpha_i
\subset \tilde{K}$ begins at $q_i$ and ends at some $q_{\mu_i}$.  Then
$\tilde{K}$ {\it admits a braid-index reducing move} if for some
$i\in\{1,\dots,n\}$ there exists an arc
$\delta$ in $P\times\{0\}$ which joins $q_i$ to $q_{\mu_i}$, such that the
closed curve $\alpha_i\cup\delta_i$ bounds a disc $\Delta$ which intersects
$\tilde{K}$ precisely in the braid strand $\alpha_i$. See sketches (a) and (d)
of Figure \ref{figure:reducible1} for an example.
The caption describes the example in terms of the standard elementary braid
generators
$\sigma_1,\dots,\sigma_{n-1}$ for ${\bf B}_n$.\bs

If one deforms $\delta$ in $P\times\{0\}$ to $\delta^\prime$, where $\delta'$
has a point on ${\bf A}$, as in sketch (d), then $\delta^\prime$ will
sweep out the disc $\Delta$ as it is moved around the axis, keeping one
endpoint fixed on ${\bf A}$ and sliding the other along $\tilde{K}$ to
$\beta(\delta^\prime)$. Thus, if $\delta$ exists, then the closed $n$-braid
$\tilde{K}$ can be replaced by a new representative of the same knot type,
$\tilde{K} - \alpha_i + \delta$, by pushing $\alpha_i$ across $\Delta$ to
$\delta$.  Moreover, after tilting $\delta$
slightly to make it transverse to the fibration, the modified representative
$\tilde{K'}$ will be an $(n-1)$-braid representative of the same knot type.
See sketch (b).  Thus the closed $n$-braid $\tilde{K}$ is {\it reducible} to a
closed
$(n-1)$-braid, and so we say that it {\it admits a reducing move} or {\it has a
reducing loop}.  If we orient the braid strands anticlockwise,  and
number them along the positive x axis  so that the one that is furthest from
the
braid axis $A$ is strand 1, then in the example $n=4, i=1$ and $\mu_i=3$.
Remark: this numbering is correct because each $P\times\{t\}$ has a
natural orientation, when we think of it as a disc in $S^3$, i.e. the
orientation determined by the orientation on ${\bf A} =
\partial (P\times\{t\})$. In Figure
\ref{figure:reducible1}(d) the z-axis and the positive x-axis determine the
plane  $P\times\{0\}$, and in sketch (d) we are seeing it from its negative
side.  \bs

Observe that the
property `$K$ admits a reducing loop' is a property of the conjugacy class of
$K$. For, the reducing arc $\alpha_i$ meets every fiber
$P\times\{t\}, \ \ t\in [0,1]$, so that we can find a reducing arc
$\alpha_i(t)$ by cutting
$\tilde{K}$ along $P\times\{t\}$ if and only if we can find a reducing arc
$\alpha_i = \alpha_i(0)$ by cutting along $P\times\{0\}$.

\begin{figure}[htpb!]
\begin{center}
\centerline{\includegraphics[scale=.9, bb=98 587 453 673]{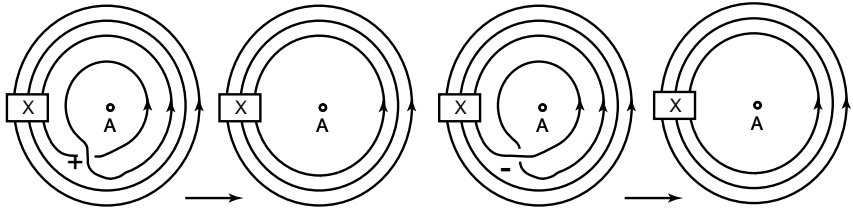}}
\caption{Positive and negative reducing loops}
\label{figure:reducible2}
\end{center}
\end{figure}

We next define the {\it sign} of a reducing loop. Assume that $\tilde{K}$
admits a reducing move at the arc $\delta$, as in Figure
\ref{figure:reducible1} (a) or \ref{figure:reducible2}. After deleting the
reducing loop the algebraic crossing number of the braid, i.e. its exponent sum
when
described by a word in the standard generators of ${\bf B}_n$,   will either
decrease or
increase by 1. If it decreases (resp. increases) we say that the reducing loop
was positive (resp. negative). The left and right examples in Figure
\ref{figure:reducible2} are positive and negative respectively. The reducing
loop in Figure \ref{figure:reducible1}(a) is positive because 
the algebraic crossing
number goes down by 1 after the reduction.\bs

\begin{figure}[htpb!]
\centerline{\includegraphics[scale=1.0, bb=162 559 386 645]{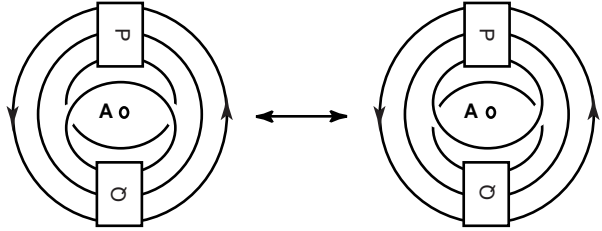}}
\caption{The exchange move} 
\label{figure:exchange1}
\end{figure}

The conjugacy class $\tilde{K}$ is
said to admit an {\it exchange move} if it has a representative $K$ which has
the special form $P\sigma_{n-1}Q\sigma_{n-1}^{-1}$, where $P$ and $Q$ only use
strands $1,\dots,n-1$.  See Figure \ref{figure:exchange1}. Note that if $K$
admits
an exchange move, then $K$ is a product of two reducible braids,
$P\sigma_{n-1}$
 and
$Q\sigma_{n-1}^{-1}$ of opposite sign. \bs

The knot or link type
${\cal K}$ is said to be {\it exchange reducible} if, up to conjugacy, any
closed braid representative can be reduced to any representative of
minimum braid index by a sequence of reducing  and exchange moves. \bs

\begin{figure}[htpb!]
\begin{center}
\centerline{\includegraphics[scale=1.0, bb=77 374 401 531]{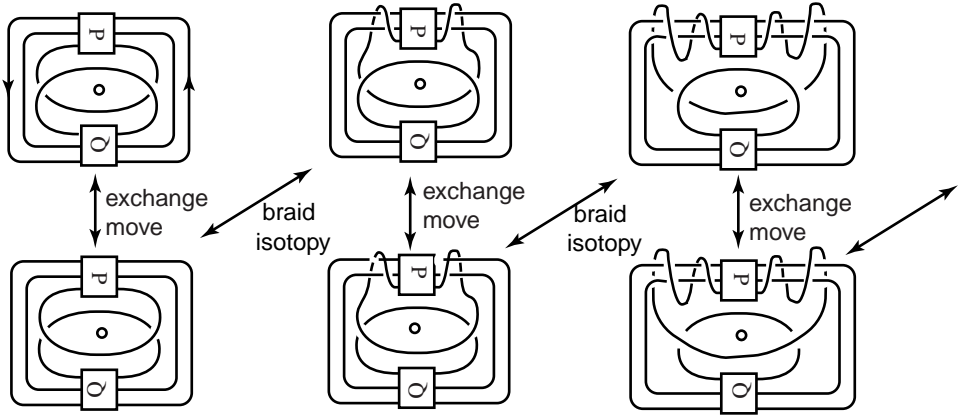}}
\caption{Exchange moves and braid isotopy can lead to infinitely many conjugacy
classes of closed braid representatives of a knot or link}
\label{figure:exchange2}
\end{center}
\end{figure}

Here are
some reasons why reducing and exchange moves are of interest in algorithmic
knot theory:
\be
\item  In \cite{BiMe-IV} it was shown that if ${\cal K}$ is a split or
composite link type, then exchange moves alone suffice to modify any closed
braid
representative
$\tilde{K}$ to a split or composite closed braid.
\item It was proved in \cite{BiMe-V} that the unlink on any number of
components is exchange-reducible.
\item  It was proved In \cite{Men} that iterated torus knots and links are
exchange-reducible.
\item Exchange moves (together with isotopy in the complement of the braid
axis) can lead to a major complication in the study of knots via closed braids,
namely exchange moves result in infinitely many conjugacy classes
of closed $n$-braid representatives of a knot or link, all related by exchange
moves.
See Figure \ref{figure:exchange2}.  Fortunately (see \cite{BiMe-VI}) it has
been
proved that if a knot or link has infinitely many conjugacy classes of n-braid
representatives then all but finitely many of them are related by
exchange moves. Since it is fairly obvious from the pictures in Figure
\ref{figure:exchange2} that there ought to be a class of `minimum complexity',
it
then becomes very important to recognize when a conjugacy class admits an
exchange
move.
\ee
In this paper we explore the
question: how can we recognize whether the conjugacy class of a given closed
braid admits a reducing move or an exchange move? \bs

The question of detecting reducing moves was posed in \cite{BLMG}. In
\cite{McCool} McCool proved the existence of an algorithm for determining
whether the conjugacy class of a braid admits a reducing move, using Garside's
solution to the word and conjugacy problems in ${\bf B_n}$.   The problems of
recognizing reducing moves and exchange moves  were considered by Fehrenbach in
his (unpublished) PhD thesis \cite{FL}, using an approach which is very
different from ours. In
\cite{Fiedler} Fiedler gives a partial invariant which shows that the conjugacy
class of a braid of braid index
$\geq 4$ is in general changed by exchange moves.  (For $n=3$ exchange moves
can be realized by conjugacy, except in the special case of composite knots).
To the best of our knowledge no algorithmic solution exists to the problem of
recognizing when a braid conjugacy class admits an exchange move.  \bs

At the same time that the connections which we just described
between braids and links were under investigation, a parallel
but quite different set of investigations was revealing new information about
matrix representations of the family of braid groups $\{{\bf B}_n, \ n\in
\ints^+\}$. For many years a central problem about the braid groups had been
the
question
of whether they were or were not linear groups.  Inspired by the work of
Thurston (see \cite{FLP}) on surface mapping class groups, it was shown in
\cite{BLMcC}, \cite{Mc} and also \cite{Iv} that braid groups shared deep
structural properties with linear groups.  However, there was more to it
than that. In the 1930's an especially interesting irreducible
representation of
${\bf B}_n$ in $GL_{n-1}(\ints[t,t^{-1}])$ had been introduced by W. Burau
\cite{Bu}, and this {\it Burau representation} of ${\bf B}_n$ was widely
regarded as a likely candidate for the sought-for faithful matrix
representation.  It was recognized by various people that the representation
had a natural interpretation via the action of ${\bf B}_n$ on the
$n$-times punctured plane $P_n$ which lifted to an action on
$H_1(\tilde{P_n},\ints)$ of the infinite cyclic covering space
$\tilde{P_n}$ of $P_n$, where the covering translation gave a module structure
to $H_1$. This fact enabled the second author of this paper to prove, in
\cite{Moo}, that for sufficiently large $n$ the Burau representation is in fact
{\it not} faithful.
The key new fact which made that proof possible was the interpretation of the
entries in the Burau matrix of a geometric braid $K$ as
recording information about intersections between the lifts of
certain arcs on $P_n$ and their images under a braid homeomorphism $\beta$ to
$\tilde{P_n}$, and that the proof of faithfulness rested on whether one could
construct  two  arcs
which have essential geometric intersections, but tricked the Burau
matrix into thinking they did not intersect.
Non-faithfulness was then proved for
$n\geq 9$ (resp. 6,5) in
\cite{Moo} (resp.
\cite{LP}, \cite{Big1}).     \bs

At this juncture there may be a historic opportunity, as
developments in both the geometry and the algebra are beginning
to bridge the gap between the two techniques. Regarding the
geometry, it is now possible to replace the Markov theory with
the more powerful Birman-Menasco theory, a relevant
portion of which was just  described.   As for the algebra,
it is now known that there is no shortage of representations of
${\bf B}_n$ in subgroups of the linear group over a
ring of Laurent polynomials (in general with more than 1 variable), moreover
they
sort themselves out into families in a way that suggests close connections with
the
Burau representation, but with one important difference: the newly discovered
representations include ones that are faithful.  We have in mind first the
important
work in Ruth Lawrence's thesis
\cite{Law}. Second, we are thinking of the key contribution from Daan Krammer,
who rediscovered the Lawrence representation (in a different setting) and used
it to prove the linearity of ${\bf B}_4$ in \cite{Kra1}. Third, we have in mind
the proof by S. Bigelow in \cite{Big} and D. Krammer in \cite {Kra2} that the
`Krammer Lawrence' or KL representation is faithful for all $n$.
Finally, we have in mind
the construction which was first suggested in \cite{BLM} and then
investigated in detail in \cite{Lo}. As will be proved below, in the simplest
case the latter construction gives a variation on the KL representation.
Like the Burau representation, it detects intersections of
certain arcs on the
$n$-times punctures disc, and this is the fact that we exploit in relation to
the
problems in algorithmic knot and link theory that we discussed above.   \bs

The main result of this paper will be the development of algebraic techniques
to
detect when the conjugacy class of a closed braid which represents a knot
admits a reducing or exchange move. See Theorems
\ref{theorem:detecting reducibility} and \ref{theorem:detecting exchange moves}
and Corollary \ref{corollary:recognizing exchange moves}.
Here is the plan of this paper. In Section \ref{section:representations of Bn}
we establish our conventions and define the particular
representations that we shall use here.
In section
\ref{section:an intersection pairing}  we define the intersection
pairing and determine its properties. In Section
\ref{section:detecting reducing moves} (resp. \ref{section:detecting exchange
moves}) we establish the relationship between the intersection pairing created
by the
action of a braid homeomorphism ${\beta}$ on the n-times punctured plane and
the
detection of reducing loops (resp. exchange moves) in the associated closed
braid.
 In Section
\ref{section:detecting intersections via the monodromy matrices},
we establish
the connection between the intersection pairing and certain blocks of zeros in
the image of $\beta$ under our matrix representation of ${\bf B}_n$.
In Section \ref{section:questions, comments and conjectures} we  explain why
our
work does not, at this time, give an algorithm for recognizing the unknot. We
discuss the open problems whose solution would make it into an algorithm. We
also
discuss several other interesting open problems which were suggested by the
work
 in
this paper.

\section{Representations of ${\bf B}_n$}
\label{section:representations of Bn}
In this section we'll review the general
construction of  `homology representations' of ${\bf B}_n$. We begin by
establishing our conventions and notation. After that we construct the
representation
which is of primary interest, first in a special case, and then in a more
general
setting. Finally, we show why the special case gives a reducible  form of the
KL representation.

\subsection{Braid homeomorphisms and geometric braids}
\label{subsection:braid homeomorphisms and geometric braids}
We will need to use two ways of looking at a closed braid: as a
geometric closed braid $\tilde{K}$ in
$\reals^3$ (with a particular choice of a half-plane $P$ along which we cut it
open
to a braid
$K$) and as a homeomorphisms ${\beta}$ of $P$, punctured at $n$ points. We now
explain our conventions for passing between $\beta$ and $K$. Let
$$P = \{(x,y,z) \in {\mathbb R}^3:y=0, x\ge 0\}$$
and let $A\subset P$ be the $z$-axis.  Then  there is a map
$F: P \times [0,1] \to {\mathbb R}^3$  defined by
$$F((x, 0, z),t) =(x cos(2\pi t), x sin (2 \pi t), z).$$
Note that when $t=0$  or 1 the restriction of $F$ to $P$ is the identity map.
\\ \bs

 Choose a set of points $Q_n = \{q_1,\dots,q_n\}$ on $P$ in the
interior of a disc $D^2\subset P$.  Let ${\beta}:P\to P$ be a {\it braid
homeomorphism}, i.e. a homeomorphism of $P$ which fixes the set $Q_n$ and is
the
identity on $P\setminus int(D)$.  Then ${\beta}$ is isotopic to the
identity as a homeomorphism of $P$ and so there is an isotopy (which we may
assume is fixed on $P\setminus int(D)$), say $H:P\times [0,1]\to P$, with
$H(p,0) = p, H(p,1) ={\beta}( p)$ for all $p\in P$. Then $H$ induces a
homeomorphism
$H':P\times [0,1]\to P\times [0,1]$ which is defined by  $H'(p,t) =
(H(p,1-t),1-t)$.  The composite map ${\bf H} = F\circ H'$ induces
identifications
 ${\bf H}(p,0) ={\bf H} (p, t)$   if   $p\in A$ and ${\bf H}(p,1)
={\bf H}({\beta}(p),0)$  for all $p\in P$.
Clearly ${\mathbb R}^3$, Euclidean space, is obtained from $P \times
[0,1]$ by these
identifications. The image of  $Q_n\times [0,1]$ under ${\bf H}$ will
be the closed geometric braid $\tilde{K}$ {determined
by} the braid homeomorphism ${\beta}$. The fact that we have chosen the
half-plane
$P\times \{0\} = P\times \{1\}$ as our reference gives us, in a natural way,
an open geometric braid associated to ${\beta}$: it is obtained
by cutting $\tilde{K}$ along the half plane $P\times \{0\}$. See the example in
Figure \ref{figure:reducible1}.   \bs

Now that we have learned how to pass between $\beta$ and $K$, we will not need
to
distinguish between them.   Thus we may pass back and forth freely between the
interpretation of our braid as a homeomorphism of the n-times punctured plane
and as a geometric braid in 3-space. We will use the symbol $\beta$ for
elements
of ${\bf B}_n$ and the symbol $\tilde{\beta}$ for the conjugacy class of
$\beta$. This is the same as the isotopy class of the
associated closed braid, where isotopy means isotopy in the complement in
$\reals^3$ of the braid axis.
\subsection{The group ${\bf B}_{1,n}$}
\label{subsection:The group B}  We will be interested in the braid
group ${\bf B}_n$ on n-strands, but in order to look at ${\bf B}_n$ as a group
of
isotopy classes of homeomorphisms of the punctured plane it will be convenient
to
regard it as a subgroup of
${\bf B}_{n+1}$. Number the strands in the latter group as $0,1,\dots,n$.   Let
${\bf B}_{1,n}\subset {\bf B}_{n+1}$ be the subgroup of braids in
${\bf B}_{n+1}$ whose associated permutation fixes the letter 0. Its
relationship to ${\bf B}_n$ is given by the  group extension
\begin{equation}
\label{short exact sequence defining B(1,n)}
 1\to {\bf F}_n \to {\bf B}_{1,n}\to {\bf B}_n \to 1,
\end{equation}
where the homomorphism ${\bf B}_{1,n}\to {\bf B}_n$ is defined
by pulling out the zero$^{th}$
braid strand.  There is a cross section which is defined by
mapping ${\bf B}_n$
to the subgroup of braids on strands
$1,\dots,n$ in
${\bf B}_{1,n}$. Therefore we may identify ${\bf B}_{1,n}$ with
${\bf F}_n\semi{\bf B}_n$. \\ \bs

We use bold-faced roman letters for elements of ${\bf F_n}$, Greek
letters for elements of ${\bf B_n}$, and upper case Roman
letters for elements of
${\bf B}_{1,n}$. Since
${\bf F}_n$ and
${\bf B}_n$ are both known  groups, and since the action
of  ${\bf B}_n$ of ${\bf F}_n$ is also known, we
conclude that ${\bf B}_{1,n}$ has a
presentation with generators
${\bf x_1,...,x_n}, \sigma_1,...,\sigma_{n-1}$, where the ${\bf x_i}'s$
generate 
the free factor ${\bf F}_n$ and the elementary braids $\sigma_j$ generate
the factor ${\bf B}_n$, satisfying the well-known braid relations:
\begin{equation}
\label{equation:braid relations}
\sigma_i\sigma_j = \sigma_j\sigma_i \ \ \ \ {\rm if} \ \ \ |i-j|\geq 2,
\ \ \ \ \ \ \ {\rm and} \ \ \ \ \ \ \ \ \sigma_i\sigma_j\sigma_i =
\sigma_j\sigma_i\sigma_j
\
\
\
\  {\rm if}
\
\ |i-j| = 2.
\end{equation}
If a word $\sigma_{\mu_1}^{\epsilon_1}\cdots\sigma_{\mu_r}^{\epsilon_r}$ 
represents an element of ${\bf B_n}$, then $\sigma_{\mu_1}^{\epsilon_1}$
corresponds to the first crossing in the oriented braid. \bs 

Our conventions for
${\bf x_1,\dots,x_n}$ are that the generator
${\bf x_i}$ represents a braid in which strand 0 begins at $q_i\times\{0\}$,
travels in front of strands $1\dots,i-1$, then around strand
$i$ with positive linking number, and then passes to $q_i\times\{1\}$
traveling in front of strands $1\dots,i-1$ again.
The geometry then determines additional relations which give the action of
${\bf B}_n$ on ${\bf F}_n$ (\ref{equation:action of Bn
on Fn}):
\begin{equation}
\label{equation:action of Bn on Fn}
\sigma_i{\bf x_j}\sigma_i^{-1}  = 
\cases {{\bf x_{i+1}} & if  $j=i$;\cr
 {\bf x_{i+1}}^{-1}{\bf x_i}{\bf x_{i+1}} & if  j=i+1\cr
{\bf x_j} & otherwise}
\end{equation}
In this way we see that the groups ${\bf B}_n$ and ${\bf F}_n$ are
both subgroups of ${\bf B}_{1,n} \subset {\bf B}_{n+1}$.  \bs

\subsection {Magnus representations and an example}
\label{subsection:Magnus representations and an example}
In this section we show how the ideas which were introduced by W. Magnus to
explain the underlying mechanism behind the Burau representation of ${\bf B_n}$
(see \cite{BLMG}) can be bootstrapped to yield new representations of
${\bf B_n}$. Later we will see that our new, augmented representation is in
fact
the Lawrence-Krammer representation of ${\bf B_n}$. \bs
 
Recall that the free calculus was used 
by Gassner in \cite{Ga} to construct an irreducible representation of the pure
braid
group
${\bf P}_{n+1}$ in the ring of
$n$-dimensional matrices whose entries are Laurent polynomials in $n$
variables. In \cite{BLMG} it was explained how similar ideas could be
used to construct representations of the subgroups of
${\bf B}_{n+1}$ which lie between ${\bf B}_{n+1}$ and  ${\bf P}_{n+1}$. The one
that
interests us here is the irreducible  representation of the in-between subgroup
${\bf B}_{1,n} = {\bf F}_n\semi{\bf B}_n$  over the ring
$\ints[t,t^{-1},q,q^{-1}]$ of Laurent polynomials in 2 variables. The variable
$t$ is to be thought of as being associated with strand 0 and
the variable $q$ as being associated to strands 1 through $n$.  Call this the
{\it Magnus} representation of ${\bf B}_{1,n}$.    Let ${\cal M}$ be the ring
of
$n+1$ by $n+1$ matrices with entries in $\ints[t,t^{-1},q,q^{-1}]$. The Magnus
representation is a representation
$$\rho:{\bf F}_n\semi{\bf B}_n \to {\cal M}.$$
It has dimension $(n+1)\times (n+1),$ but as is well known it's is reducible to
$n\times n$. For our purposes it will be most convenient to use the unreduced
form, and we shall do so. \\ \bs

It is shown in  \cite{BLMG} how to construct the images of the generators of
${\bf F}_n\semi{\bf B}_n$ under $\rho$. Working out an example, we obtain
the following matrices for the generators of $\rho({\bf F}_n\semi{\bf B}_n)$ in
the
case $n=4$:
$$ \sigma_1\mapsto\left(\matrix 
{1&0&0&0&0\cr 0&0&q&0&0\cr 0&1&1-q&0&0\cr
0&0&0&1&0\cr 0&0&0&0&1}\right), \ \ \ \
\sigma_2\mapsto
\left(\matrix{1&0&0&0&0\cr 0&1&0&0&0\cr 0&0&0&q&0\cr 0&0&1&1-q&0\cr
0&0&0&0&1}\right)$$
$$\sigma_3\mapsto\left(\matrix{1&0&0&0&0\cr 0&1&0&0&0\cr 0&0&1&0&0\cr
0&0&0&0&q\cr
0&0&0&1&1-q}\right)$$

and \\

\

$$ x_1 \mapsto \left(\matrix{
q&-q(1-q)&0&0&0\cr
1-t&1+tq-q&0&0&0\cr
0&0&1&0&0\cr
0&0&0&1&0\cr
0&0&0&0&1}\right) ,$$

$$ x_2 \mapsto \left(\matrix{
q&-(1-q)^2&-q(1-q)&0&0\cr
0&1&0&0&0\cr
1-t&(1-t)(1-q)&1-q+tq&0&0\cr
0&0&0&1&0\cr 0&0&0&0&1}\right), $$

$$x_3 \mapsto \left(\matrix
{q&-(1-q)^2&-(1-q)^2&-q(1-q)&0\cr
0&1&0&0&0
\cr 0&0&1&0&0\cr
1-t&(1-t)(1-q)&(1-t)(1-q)&1-q+tq&0\cr
0&0&0&0&1}\right),
$$

$$ x_4 \to \left(\matrix{
q&-(1-q)^2&-(1-q)^2&-(1-q)^2&-q(1-q)\cr
0&1&0&0&0\cr
0&0&1&0&0\cr
0&0&0&1&0\cr
1-t&(1-t)(1-q)&(1-t)(1-q)&(1-t)(1-q)&1-q+tq}\right). $$

\

Let ${\cal D}(z)$ denote a diagonal matrix whose diagonal entrises are $z$. We
next define the homomorphism
\begin{equation}
\tau:{\bf F}_n\semi{\bf B}_n \to {\cal M} \ \ \ \ {\rm by} \ \ \ \
\tau(G) = {\cal D}(q^{{\rm deg}(G)})\rho(G),
\end{equation}
where $G\in {\bf B}_{1,n}$ and deg$(G)$ denotes the degree of $G$, i.e. its
exponent sum in the basis elements ${\bf x_i}$ of the free group ${\bf F_n}$.
In
particular,
$\tau(\sigma_i) = \rho(\sigma_i)$,   whereas
$\tau({\bf x_j}) = {\cal D}(q)\rho({\bf x_j})$. \bs

Using the representation $\tau:{\bf B}_{1,n}\to {\cal M}$, we are finally
ready to construct (in our special case) the representation $\tau^+$ of ${\bf
B}_n$
which is of primary interest to us in this paper.
It is essentially identical to Lawrence's version (see [18]) of the now-famous
Krammer-Lawrence representation.  Let $P_n = P - Q_n$ denote the $n$-times
punctured plane. View $F_n$ as the fundamental group of $P_n$ based at some
point
$x \in A = \partial P_n$, and use the representation $\tau$
(restricted to the first factor) to define a local system ${\cal L}$ on $P_n.$
We already gave the action of ${\bf B_n}$ on ${\bf F_n}$ in
(\ref{equation:action of Bn on Fn}). It will be helpful to think
of a  local system as a covering space which has the structure of a fiber
bundle, where the fibers are copies of the Magnus reprsentation
$V={\mathbb Z}[t,t^{-1},q,q^{-1}]^n.$ The homology $H_1(P_n,x,{\cal L})$
is sometimes called $H_1(P_n,x,V),$ homology with coefficients
in $V$ viewed as a module over the
fundamental group of $P_n$. It is a direct sum $V\oplus V \oplus ... \oplus V$
($n$ factors) and we denote by $e_i:V \to H_1(P_n,x,{\cal L})$ the inclusion on
the $i$'th factor.
For $\alpha \in V$ the identification is such that
$e_i(\alpha)$ is the loop $x_i$ lifted to the cover, in such a way that the
basepoint $x$
at the beginning lifts to the element $\alpha \in V.$ \bs

We use the action (\ref{equation:action of Bn on Fn}) to compute the map
induced by the elementary braid generators $\sigma_i$ on homology:
$$\sigma_i:H_1(P_n, x, {\cal L}) \to H_1(P_n, x, (\sigma_i)_*({\cal L})).$$
 Since
it is enough to determine the homology classes $\sigma_i(e_j)(\alpha)$
for $\alpha \in V$ it is enough to calcluate the composites $\sigma_i e_j.$
Identifying both source and target of the map
with $V\oplus V \oplus ... \oplus V$ in the manner described
above,  we find:

\begin{equation}
\label{equation:action of Bn on homology}
\sigma_i e_j   = \cases{ e_{j+1}(\tau({\bf x}_{j+1}))& if $j=i$ \cr
 e_i(\tau({\bf x}_i)) + e_{i+1}(1-\tau({\bf x}_{i+1}))& if $j=i+1$\cr
e_j & otherwise}
\end{equation}
To understand (\ref{equation:action of Bn on homology}), we think of
$e_{i+1}(\alpha) $ as the loop
$x_{i+1}$ lifted to the cover, in such a way that the base point $x$ at the
beginning lifts to the element $\alpha\in V.$ When we apply $\sigma_i$ we get a
loop that traces
$x_{i+1}$ anticlockwise, then $x_i$ anticlockwise, and finally $x_{i+1}$
clockwise.
However, the initial points of these loops lift to three different
sheets in the cover, i.e. to three different elements of elements of
$V$. These are the coefficients of  $e_{i+1}, e_i, $ and $
-e_{i+1}$, taking account of the fact that the action of a loop $w$ in $P_n$
on $(\sigma_i)_*({\cal L})$ is the action of $\sigma_i(w)$ on ${\cal L}.$ \bs

Note that in future sections we will make the
following identification. Since
${\cal M}=End(V)$ is the endomorphism ring of
$V,$ The map
$e_i:V\to H_1(P_n,x,V)$ can be viewed as an element of $H_1(P_n,x, {\cal M})$
where the homology is taken with respect to the left action of $F_n$ on ${\cal
M}$. This homology module is a free right ${\cal M}$ module of rank $n$ and the
$e_i$ are a basis for this module. The formula displayed above describes
the action of $\sigma_i$  with respect to this basis. However we are not
quite done. 
The resulting matrix with ${\cal M}$ entries
does not define a representation because our identification
of the source and target both with $V \oplus V \oplus ... \oplus V$
was not natural.  \bs  

For any braid $\beta$ there is a natural map  
$\beta_*{\cal L}\to {\cal L}$
which is given on the fiber $V$ over $x$ by multiplication
by $\tau(\beta),$ and it extends to a map of the whole local
system because of the relations (3) in $ {\bf B}_{1,n}.$ 
Now let $\alpha$ be any other braid. Applying $\alpha_*$ we obtain
a second map
of local systems $\alpha_*\beta_*{\cal L}\to \alpha_*{\cal L}$
Because homology is a bifunctor we then obtain a commutative square
$$\matrix{  
H_1(P_n,x,{\cal L})&\to&H_1(P_n,x,\beta_*{\cal L})&\to&
H_1(P_n,x,\alpha_* \beta_*{\cal
L}) \cr
                  &   & \downarrow & & \downarrow \cr
                  &   &H_1(P_n,x,{\cal L})&\to&H_1(P_n,x,\alpha_*{\cal L})\cr
                   &   &                    & &\downarrow\cr
                  &   &                    & &H_1(P_n,x,{\cal L})  }.$$
Each module in the diagram is a direct sum of $n$ copies of $V$
and so each map may be represented by an $n$ by $n$ matrix with entries in
${\cal M}. $ The leftmost vertical map  in the square is  given by a diagonal
matrix ${\cal D}(\tau(\beta))$ with each entry equal to
$\tau(\beta).$ The horizontal maps are the maps on homology
given on generators by formula (5). The fact that the square
commutes implies an equality of two maps from the upper left corner to the
lower right corner
${\cal D}(\tau(\alpha))\alpha {\cal D}(\tau(\beta))\beta
= {\cal D}(\tau(\alpha\beta))\alpha\beta,$ so the assignment of each
braid generator $ \sigma_i$ to the matrix ${\cal D}(\sigma_i)\sigma_i$
describes a braid group representation.
Thus, if ${\cal D}(\tau(\sigma_i)), \ i=1,2,3$ denotes the block diagonal
matrix (with each block a $5 \times 5$ matrix) whose diagonal entries are the
blocks
$\tau(\sigma_i)$, defined earlier, the `augmented' representation
 $\tau^+$ of ${\bf B}_4$ is given on the elementary braid generators by:

$$ \tau^+(\sigma_1) = {\cal D}(\tau(\sigma_1))\left(\matrix
{0&\tau(x_1)&0&0\cr
1&1-\tau(x_2)&0&0\cr  0&0&1&0\cr
0&0&0&1}
\right) =$$
$$\left(\matrix
{0&\tau(\sigma_1)\tau(x_1)&0&0\cr
\tau(\sigma_1)&\tau(\sigma_1)(1-\tau(x_2))&0&0\cr  0&0&\tau(\sigma_1)&0\cr
0&0&0&\tau(\sigma_1)}
\right)
$$

$$ \tau^+(\sigma_2 )={\cal D}(\tau(\sigma_2))
\left(\matrix
{1&0&0&0\cr 0&0&\tau(x_2)&0\cr
0&1&1-\tau(x_3)&0\cr
0&0&0&1}
\right)
=$$ 
$$\left(\matrix
{\tau(\sigma_2)&0&0&0\cr 0&0&\tau(\sigma_2)\tau(x_2)&0\cr
0&\tau(\sigma_2)&\tau(\sigma_2)(1-\tau(x_3))&0\cr
0&0&0&\tau(\sigma_2)}
\right)
$$

$$ \tau^+(\sigma_3)= {\cal D}(\tau(\sigma_3))\left(\matrix
{1&0&0&0\cr
0&1&0&0\cr
0&0&0&\tau(x_3)\cr
0&0&1&1-\tau(x_4)} \right)
=$$
$$ \left(\matrix
{\tau(\sigma_3)&0&0&0\cr
0&\tau(\sigma_3)&0&0\cr
0&0&0&\tau(\sigma_3)\tau(x_3)\cr 
0&0&\tau(\sigma_3)&\tau(\sigma_3)(1-\tau(x_4))} \right).
$$
The rank of this representation is  $4\times 5 = 20$ (for general $n$ it is
$n(n+1)).$  If we had used the reduced Magnus matrices the rank would have been
$n^2.$  \bs

\subsection{The general case}
\label{subsection:the general case} 
Suppose we are given a representation $\tau:{\bf B}_{1,n}\to
GL(V)$. The module $H_1(P_n,x,{\cal L})$ has a
${\bf B}_n$ action, coming from the  fibration
$$P_n \to X_{1,n}\to X_n$$
where $X_{1,n}$ is the space of $n$ distinct points in $P$ with one marked
point, and $X_n$ is the space of $n$ distinct points in $P.$ The second map
is forgetting the marked point. The fiber over a subset $Q_n \subset P$
is $P_n.$  The  action of ${\bf B}_{1,n}=\pi_1(X_{1,n})$ on $V$ is
taken as the monodromy defining a local system on $X_{1,n}$
and as usual the fundamental group of
the base acts on the homology of the fiber over $Q_n.$
Note that the group $H_1(P_n,x, {\cal L})$ is homology relative
to the basepoint.  
This is a direct sum of $n$ copies of $V$. The
action of a braid $\beta$ is the composite of the homology action
$H_1(P_n, x, {\cal L}) \to H_1(P_n,x,\beta_*{\cal L})$   with
the map induced on homology by
$\beta_*{\cal L}\to {\cal L}$,  which  multiplies the fiber $V$ over $x$ by
$\tau(\beta).$  Thus, given a representation $\tau:{\bf B}_{1,n}\to GL(V)$ one
obtains immediately  a new representation $\tau^+:{\bf B}_n\to GL(V\oplus ...
\oplus
V)$.  

\subsection{The Krammer-Lawrence representation}
\label{subsection:the Krammer-Lawrence representation}
We now show that our $n(n+1)$ dimensional representation contains the one
Krammer and Bigelow have proven faithful (with our $`t$' playing the
role of Krammer's $t^2$). Although we've displayed the unreduced Magnus
matrices
above, for the moment let us take $\rho:{\bf B}_{1,n}\to {\cal M}$ to be
the size $n$ reduced matrices.  Let $X_n$ be the set of $n$ element subsets
$Q_n\subset P,$ let  
 $X_{1,n}$ be the set of subsets $\{q_0,...,q_n\}\subset P$
with $q_0$ marked, and let $X_{2,n}$ be the set of subsets
of order $n+2$ of $P$ with two separately
marked points. Recall $Q_n=\{q_1,...,q_n\}$ and $P_n=P-Q_n.$
Let $SP$ be the space of pairs of distinct, separately marked points
in $P_n.$
Let ${\cal L}$ be the local system of $X_{2,n}$ with
fiber ${\mathbb Z}[t,t^{-1},q,q^{-1}]$
and monodromy
matrices as in Krammer's and Bigelow's papers, so that the power of $q$ is the
winding number of either marked point around $q_1,...,q_n$ and
the power of $t$ is the winding number between the marked points.
We have the diagram of fibrations
$$\matrix{(P-\{q_0,..,q_n\})&\to& SP & \to &P- \{q_1,..,q_n\} \cr
           \downarrow && i\downarrow && k\downarrow\cr
           (P-\{q_0,...,q_n\})&\buildrel j\over
\rightarrow&X_{2,n}&\to&X_{1,n}\cr
&&\downarrow&&\downarrow\cr
&&X_n&\to&X_n\cr}$$
The representation of $B_n$ defined by $\tau^+$ maps onto $V$ and the kernel
is
a copy of the unreduced homology
$H_1(P_n, k^* {\cal H})$ where ${\cal H}$ is the
local system on  $X_{1,n}$ corresponding to the representation
$\tau:{\bf B}_{1,n}\to {\cal M}.$  We wish to compare this with
the Krammer Lawrence representation of $B_n.$
From the spectral sequence of a fibration we have an isomorphism
of  representations of  the fundamental group $B_n$ of $X_n$.
\begin{equation} \label{equation:isomorphism of local systems}
 H_2(SP, i^*{\cal L})=
 H_1(P_n, k^*{\cal H}_1(P-\{q_0,...,q_n\},j^*{\cal L}))
\end{equation}
The module on the left is the Krammer-Lawrence
representation (but with two marked points instead of two indistinguishable
points)
while module on the
right is
 the homology $H_1(P_n,k^* {\cal H}')$ with coefficients
in the local system $
{\cal H}'={\cal H}_1(P-\{q_0,...,q_n\},j^*{\cal L}).$
To prove $H_1(P_n,k^*{\cal H})$ is isomorphic to $H_1(P_n,k^*{\cal H}')$
it now remains only to
identify the two local systems  ${\cal H}$  and  ${\cal H}'$  on $X_{1,n}$.
For $G\in {\bf B}_{1,n}$ the monodromy of $G$  on the first local system is
$\tau(G).$ The monodromy on the second is
the composite of the homology map
$H_1(P-\{q_0,...,q_n\}, j^*{\cal L})\to H_1(P-\{q_0,...,q_n\},g_*j^*{\cal
L}),$ represented by the reduced Magnus matrix $\rho(G)$,
followed by the map on homology induced by $\beta_*j^*{\cal L}\to j^* {\cal
L}$, which is multiplication by $q^{{\rm deg}(G)}$. The composite
is $\tau(G),$  as needed.
To get the Krammer-Lawrence representation of two
{\it unmarked points}  it is necessary to adjoin
$\sqrt{t}$ to our coefficient ring (remember our
 $t$ is Krammer's $t^2$).
The local system ${\cal L}$ is now the pullback of a local system
on $X_{2,n}/C_2$ where the cyclic group $C_2$ acts freely by
interchanging the markings on the two points. 
We obtain 
a $C_2$ action on $H_2(SP)$ and  since the Krammer-Lawrence representation
has no 2-torsion so we may identify it with the
$C_2$ invariant submodule, which
is therefore faithful, from which $\tau^+$ is as well. \bs

Note too that the $\tau^+$ are all faithful (because the Lawrence-Krammer
representation has been proved to be faithful). On the other hand, the
$\tau$ are not all faithful (because the $\tau(\sigma_i)$ act on
the last
$n$ basis vectors exactly as the Burau matrices do. )  These two facts  will be
important in Theorem \ref{theorem:intersection pairing} below. \bs

\section{The intersection pairing}
\label{section:an intersection pairing}
For the rest of this section,  ${\cal M}$ will denote an arbitrary ring,
depending on 
 $n$, together with  a homomorphism
$$\tau:{\bf B}_{1,n}\to {\cal M}$$
to the units of ${\cal M}$.  We further assume that  the
rings ${\cal M}$ are appropriately nested as $n$ increases.

\begin{figure}[htpb!]
\begin{center}
\centerline{\includegraphics[scale=0.9, bb=131 557 488 663]{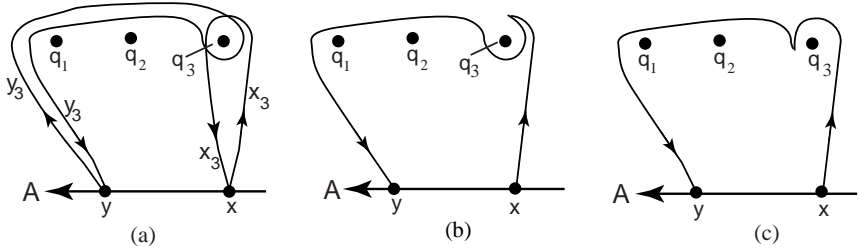}}
\caption{The intersection form}
\label{figure:intform}
\end{center}
\end{figure}

In the theorem which is stated below we will use the representation $\tau:{\bf
B}_{1,n}\to  {\cal M}$ to allow us to multiply elements of ${\cal M}$ on the
left or
on the right by elements of ${\bf B}_n$ or ${\bf F}_n.$    In order to define
fundamental classes of arcs in $P_n$  we work with homology with coefficients
in
${\cal M}.$ Choose points $x,y$ on the oriented braid axis $A.$
In the orientation of $A$ as the 
boundary of $P$ we choose $y$  next to $x$ (as in
Figure \ref{figure:intform}), but with $y$ displaced in the positive direction
compared to $x$  using
the orientation on ${\bf A}$.  Note that the axis ${\bf A}$ appears
to go from right 
to left (the clockwise direction) in the figure. This means if we use
the orientation on ${\bf A}$ to orient $P_n$ the figure shows the negative
side of the plane $P_n.$ This corresponds also to the negative
side of the plane $P_n\times\{0\}$ in figure \ref{figure:reducible1}(d),
such that the ${\bf A}$-axis in figure \ref{figure:intform}
corresponds to the  upward-pointing $z$-axis in figure \ref{figure:reducible1}
(d). In that figure the basepoints $x$ and $y$ would appear with
$y$ above $x$ on the $z$-axis. \bs

  Because the fundamental group of $P_n$ maps to the
units of the ring ${\cal M}$  it acts on any right or left
${\cal M}$ module. We denote by ${\cal M}_{left}$ the free
rank one ${\cal M}$ module upon which ${\cal M}$ acts on the
left by right ${\cal M}$ module homomorphisms,
and we denote by ${\cal M}_{right}$ the free rank one
module upon which ${\cal M}$ acts on the right by left ${\cal M}$ module
homomorphisms. Thus we consider two separate homology modules:
\begin{equation}
\label{two modules}
H_1(P_n,y,{\cal M}_{{\rm right}}) \ \ \ \ \ \ \ {\rm and}  \ \ \ \
 H_1(P_n,x,{\cal M}_{{\rm left}})
\end{equation}
The punctured plane is homotopy-equivalent
to a wedge of loops based at $x$ or $y$ which encircle and separate the
$n$ points, and since we are taking homology relative
to a point,  the module
$H_1(P_n,x,{\cal M}_{{\rm left}})$  is  equal to the cellular 1-chains of this
wedge
of circles.  This is a right ${\cal M}$ module which is free of rank $n.$ It is
a two-sided ${\bf B}_n$ module, with the right action coming from
the fact that ${\bf B}_n$ maps to ${\cal M}$ and ${\cal M}$ acts on the right,
the left action being the `interesting' homology action, which commutes with
the
right ${\cal M}$ action.
We can also form $H_1(P_n,y,{\cal M}_{{\rm right}})$ in the analogous
way. It is a two-sided ${\bf B}_n$ action with the
right action being the `interesting' one.
(Look ahead to formulas
(\ref{equation:the interesting left action on H(.., x,..)}) and
(\ref{equation:the less interesting right action on H(.., x,..)}), also
(\ref{equation:the interesting right action on H(.., y,..)}) and
(\ref{equation:the less interesting left action on H(.., y,..)}),
where the actions are worked out explicitly).
If $\gamma$ is a loop in $P_n$ based at $x$ we have
a fundamental
 homology class $[\gamma]_x\in H_1(P_n,x,{\cal M}_{{\rm left}})$ and if
$\delta$ is based at $y$ we have $[\delta]_y\in H_1(P_n,y,{\cal M}_{{\rm
right}}).$

Note that our choice of notation implies that for  
$v\in H_1(P_n,x,{\cal M}_{left}) $ and braids $\alpha, \beta \in {\bf B_n}$ the
symbol
$$ \alpha v \beta$$
means the same thing as $\tau^+(\alpha)v\tau(\beta)$ because the
$\tau^+$ representation defines the left action while the $\tau$ representation
defines the right action. Also if ${\bf x} \in {\bf F_n}$ the symbol
$$ v{\bf x}$$
will always mean $v\tau({\bf x}).$  On the other
hand if $w \in H_1(P_n,  y, {\cal M}_{right})$ then in the symbol
$$\alpha w \beta$$
it is the left action that is defined via $\tau(\alpha)$ and the right action
is one which we have not named. It will turn out to be the transform
of $\tau^+(\beta^{-1})$ via a certain correspondence $\star$ but let us
not worry about that yet. And  if ${\bf x}\in {\bf F_n}$ the symbol
$$ {\bf x}w$$ will always mean $\tau({\bf x})w.$  \bs

As before, the free group ${\bf F_n}$ will be viewed as the fundamental group
of
$P_n.$ An element ${\bf w}\in {\bf F_n}$
can be represented by a loop $w$  based at either point,
as we shall identify the fundamental group of $P_n$ based at either
$x$ or $y$ via a path in the axis $A$ connecting $x$ with $y.$ If we
choose $w$ to be a loop based at $x$ or at $y$ then we obtain  
two separate homology classes $[w]_x$ and $[w]_y.$   Note that these
belong to totally different homology groups.  We may 
omit the subscript when there is no ambiguity. The elements ${\bf x_1,...,x_n}$
will be a basis of $F_n$ represented by arcs $x_1,\dots,x_n$ which begin and
end
at $x$,  encircling the $q_i$
in such a way that the product ${\bf x_1\dots x_n}$ is represented by a
simple arc $x_1\dots x_n$ which encircles all the $q_i.$  
It is worthwhile to be quite  careful about the orientations of the
$x_i.$ In all our figures, in which we display the negative side of
$P_n$ with the axis $A$ at the bottom of the page, going from right to
left, the
$x_i$ are oriented so they appear to encircle $q_i$ in the counter-clockwise
direction.
We use the `composition convention' for composing
paths  so that in the path $x_1x_2$ the point traverses first $x_2$
and later $x_1.$  

It is convenient to
introduce a second basis ${\bf y_1,...,y_n}$ of $F_n$ where
${\bf y_i=x_1x_2...x_{i-1}x_i^{-1}x_{i-1}^{-1}...x_1^{-1}}.$ We shall usually
represent the ${\bf y_i}$ by arcs $y_i$ based at $y.$ Note there is an
involution of ${\bf F_n}$ which interchanges the roles of the ${\bf x_i}$ and
the ${\bf y_i}$.
\\ \bigskip
In the statement below,
a `loop based at $x$ (resp $y$)' will
mean any loop encircling just a single
$q_i$ which is oriented anticlockwise (resp. clockwise)
when viewed on the negative side of the plane $P$.
A `simple homology class' will mean the class of
 a  loop
 based at
$x$
(respectively $y$) which has the additional property that
it has no self-intersections. \bs

Assume that we have been given a
nested series of rings ${\cal M}$ (subscripts omitted)
and representations $\tau:{\bf B}_{1,n}\to {\cal M}$
which satisfy the hypotheses:
\begin{itemize}
\item  All
$\tau^+:{\bf B}_n\to GL_n({\cal M})$ are
faithful,
\item  The representations $\tau = \tau_n:B_{1,n}\to {\cal M}$  are
not all faithful, i.e. for some $n$ there is a non-trivial element in the
kernel of $\tau_n$,
\item  There is an anti-involution $*$ on each ${\cal M}$ such
that $\tau(G)^*=\tau(G^{-1})$ for $G\in {\bf B}_{1,n}$,
\item
The $\tau({\bf x_i})-1$ are non-zero-divisors in
${\cal M}$ which do not generate the unit
ideal in the subring of $\cal{M}$ generated by the $\tau({\bf x_i})$
and $\tau({\bf x_i})^{-1}.$
\end{itemize}

\begin{theorem}
\label{theorem:intersection pairing}

Under the hypotheses just stated, there is a generically non-degenerate 
${\cal M}$
bilinear intersection pairing
\begin{equation}
\langle\ ,\ \rangle:H_1(P_n,y, {\cal M}_{{\rm right}}\rangle\times
H_1(P_n,x,{\cal M}_{{\rm left}})\to {\cal M}
\end{equation}
which has the following properties:
\begin{enumerate}
\item [{\rm (a)}] If ${\beta}\in {\bf B}_n, \ v 
\in H_1(P_n,x,{\cal M}_{\rm left})$ and 
$w \in H_1(P_n,y,{\cal M}_{\rm right}) $ then 
$\langle{\beta} v {\beta}^{-1},{\beta} w {\beta}^{-1}\rangle ={\beta} \langle
v,w\rangle \beta^{-1}$. \\
\medskip 
If in addition $u,s \in {\cal M}$ then $\langle uv,ws\rangle =u\langle
v,w\rangle s$.

\item [{\rm (b)}] If $\beta\in {\bf B}_n, \ v \in H_1(P_n,x,{\cal M}_{\rm left})$ and $w \in H_1(P_n,y,{\cal M}_{\rm right}) $
then \ \  $\langle
v\beta,w\rangle =\langle v,\beta w\rangle.$

\item [{\rm (c)}] If $\gamma,\delta$ are loops based at $x,y$ respectively and
${\beta}\in {\bf B}_n$, then: \ \ \
$$[{\beta} \gamma]_x={\beta}[\gamma]_x{\beta}^{-1} \in H_1(P_n,x,
{\cal M}_{{\rm left}}) \ \ \ \  {\rm and} $$
$$[{\beta}\delta]_y={\beta}[\delta]_y{\beta}^{-1}
\in H_1(P_n,y,{\cal M}_{{\rm right}}).$$

\item [{\rm (d)}] There is a semilinear 
correspondence $\star$ between $H_1(P_n,y, {\cal M}_{{\rm right}})$ 
and $H_1(P_n,x,{\cal M}_{{\rm left}}).$
It takes simple homology elements to simple homology elements
and is compatible with the anti-involution
$*$ on ${\cal M}$ and with the $B_n$ action, so that for 
$r\in {\cal M}$, $v\in H_1(P_n, x, {\cal M}_{\rm left})$ 
$w\in H_1(P_n,y,{\cal M}_{\rm right})$ and $\beta \in B_n$ we have
$$( \beta vr)^\star = r^* v^\star \beta^{-1}\ \  {\rm and} $$
$$(rw\beta)^\star=\beta^{-1}w^\star r^*$$
Also for ${\bf w}\in {\bf F}_n$ if $w_x$ and $w_y$ are loops
based at $x$ and $y$ respectively which represent $w$ then
$[w_x]_x^\star=[w_y^{-1}]_y, \ \ \ \ \ {\rm and} \ \ \ \ \
 [w_y]_y^\star=[w_x^{-1}]_x$.
\item [{\rm (e)}] Let ${\bf w}
\in {\bf F}_n$. Let $w$ be a loop based at $x$ which
represents ${\bf w}.$  Let $v=[w]_x.$
Then a necessary condition for the existence of a simple arc representing
${\bf w}$ is:
$$\langle v^*,v\rangle=\tau({\bf w})-1.$$
where $1$ denotes the identity in ${\cal M}$.
\item [{\rm (f)}] If $\gamma,\delta$ are simple
loops based at $x$ and $y$ respectively then
$\langle[\delta]_y,[\gamma]_x\rangle
=0$  if and only if $\gamma$ and $\delta$ can be
homotoped, keeping their endpoints fixed, so that they do not intersect. That
is,
the pairing is effective.
\end{enumerate}
\end{theorem}

\pf  We begin our proof by constructing the intersection pairing, first via an
example which will allow us to understand both the geometry and the algebra,
and then algebraically, in full generality.
\bs

{\bf Construction of the pairing:}
We have not yet constructed the pairing, but it may be helpful to the
reader to have an example to keep in mind before we do so.  Consult Figure
\ref{figure:intform}(a).  We wish to compute
$\langle [y_3]_y,[x_3]_x \rangle$.
There are two points of intersection between $x_3$ and
$y_3$.  Travel along $x_3$, with its given
orientation,  to either intersection point and then switch to
$y_3$ to avoid the intersection, preserving orientation. This determines
an oriented arc from $x$ to $y$ for each intersection point. The arcs in the
example are shown in sketches (b) and (c). They are given by the elements
${\bf x_1x_2x_3}$ and
${\bf x_1x_2}$ of ${\bf F_n}$, where we read words from right to left
(in other words we compose paths according to the conventions of function
composition).
The
intersection is counted as positive (resp. negative) if it is necessary to make
a
right (resp. left) turn at the intersection point to avoid the intersection and
preserve orientation. In the example
$\langle [y_3]_y,[x_3]_x \rangle = \tau({\bf x_1x_2x_3}) - \tau({\bf x_1x_2})$.
Note that the intersection pairing has its values in ${\cal M}$.  \bs

We pass to the general case.  Using the basis ${\bf x_1,\dots,x_n}$ for the
free
group
${\bf F}_n$, we pass to the group ring ${\mathbb Z}{\bf F}_n$. Let $I$ be the
augmentation ideal.

As we have already observed, the homology $H_1(P_n,x,{\cal M}_{left})$
is equal to the first cellular chain module of a wedge of $n$
circles with coefficients in ${\cal M}_{left}.$ This is a free right
${\cal M}$ module of rank $n.$  The tensor product
 $I\otimes_{\bf F_{\rm n}}{\cal M}$ 
is also a free right ${\cal M}$ module
of rank $n.$ The isomorphism of free right ${\cal M}$ modules  sending the
 basis element $e_i$ of section 2.3  to
the element $({\bf x}_i-1)\otimes 1 \in I \otimes {\cal M}$
is particularly convenient
because then the action action of $B_n$ on the left
described in section 2.3 corresponds to the action  in which
an element $\iota \otimes m \in I\otimes {\cal M}$ is sent to
$$\beta( \iota \otimes m ) = \beta \iota \beta^{-1}\otimes \beta m $$
in which the conjugate $\beta \iota \beta^{-1}$ is calculated in the
(group ring of the) larger
group ${\bf B}_{1,n}.$ From now on we may as well identify
$e_i$ with $(x_i - 1)\otimes 1.$  
The homology class of an arc
$\gamma$ in $P_n$ based at $x$ is determined by its class  ${\bf g}\in {\bf
F}_n$, by the derivation
\begin{equation}
\label{equation:derivation at x}
  d:  {\bf F}_n \to \oplus e_i{\cal M}, \ \ \ \   {\bf x_i} \mapsto   e_i, \ \
\
\  d({\bf xy})=d({\bf x}){\bf y} + d({\bf y}). 
\end{equation}
The latter rule holds for all ${\bf x},{\bf y} \in {\bf F_n}$ and is the
Leibniz rule for a derivation (a Fox derivative). Note that the right
multiplication by ${\bf y}$ is by definition the right multiplication
using $\tau$ so the formula could more rigorously be written
$d({\bf xy})= d({\bf x})\tau({\bf y}) + d({\bf y}).$ Note also that
the rule implies $d({\bf x}^{-1})=-d({\bf x}){\bf x}^{-1}.$   \bs

Applying this Fox derivative to 
 the element ${\bf w}\in{\bf F}_n$ which is the homotopy class of
the loop $\gamma$   we have
\begin{equation}
[\gamma]_x = d({\bf w})\in \oplus_i e_i{\cal M}=H_1(P_n,x,{\cal M}_{{\rm
left}}).
\end{equation}
For example,
$$[x_2x_4x_2^{-1}]_x = d({\bf x_2x_4x_2^{-1}})
=d({\bf x_2 }){\bf x_4x_2^{-1}}+d({\bf x_4}){\bf x_2^{-1}}+d({\bf x_2^{-1}})$$
$$= e_2 ({\bf x_4}-1){\bf x_2^{-1}}  + e_4 {\bf x_2^{-1}}$$

This map $d$
is equivariant under the ${\bf B}_n$ action by conjugation on ${\bf F}_n$
\begin{equation}
\label{equation:property (c)}
 d({\beta} f {\beta}^{-1})={\beta}d(f){\beta}^{-1}
\end{equation}
Note that in this formula the conjugate $\beta f \beta^{-1}$ can
be calculated in $F_n$ using the equation \ref{equation:action of
Bn on Fn} while the right side of the equation can be written
as $\tau^+(\beta)d(f)\tau(\beta)^{-1}.$ \bs

The module $H_1(P_n,x,{\cal M}_{{\rm left}})$ has an ${\cal
M}$-linear  ${\bf B}_n$ action on the left. This action is the same as the
action
of ${\bf B}_n$ on
$I\otimes_{{\bf F}_n} {\cal M}$,  given by the formula ${\beta}(\iota\otimes r)
={\beta}\iota\beta^{-1}\otimes {\beta}r,$ where $\iota\in I, \ r\in {\cal M}$.
In terms of the basis $e_1,...,e_n$ the left  ${\cal M}$-linear action is given
explicitly by:
\begin{equation}
\label{equation:the interesting left action on H(.., x,..)}
         \sigma_i  e_j r   =\cases
{e_j\tau(\sigma_i) r& if $j\not= i,i+1$\cr
e_{i+1}\tau( \sigma_i) r & if $ j=i$\cr
e_i \tau(\sigma_i {\bf  x_i}) r + e_{i+1}\tau(\sigma_i)(1-\tau({\bf x_{i+1}}))
r & if $j=i+1$}
\end{equation}
while the less-interesting right $B_n$ action, which is semilinear with
respect to the $B_n$ action by conjugation on ${\cal M}$, is given
\begin{equation}
\label{equation:the less interesting right action on H(.., x,..)}
e_j  r \cdot \sigma_i =  e_j r  \tau(\sigma_i).
\end{equation}

Now we repeat everything, this time taking
homology with respect to the right action
of ${\bf B}_{1,n}$ on ${\cal M},$ using our basepoint $y$.  See Figure
\ref{figure:intform}(a). In order to describe this  module in a way which
will later make our intersection pairing simple,
we use the involution  on ${\bf F}_n$ which is given by:
$$ {\bf x_i}\mapsto {\bf y_i} =
{\bf x_1 x_2 ... x_i^{-1} x_{i-1}^{-1} ... x_1^{-1}}$$
We have that ${\bf B}_n$ acts ${\cal M}$-linearly on the right on the  left
${\cal M}$-module
$H_1(P_n,y,{\cal M}_{{\rm right}}) =
 \bigoplus_{i+1}^{n}{\cal M}f_i$  where we now take as our basis
elements the ${ f_i}=1\otimes ({\bf y_i}-1) \in {\cal M}\otimes_{\bf F_n} I.$
The right action
of ${\bf B}_n$ is easiest described by giving the action of the inverse
of the braid generators:
\begin{equation}
\label{equation:the interesting right action on H(.., y,..)}
    r f_j \sigma_i^{-1} =\cases
{ r\tau(\sigma_i^{-1}) f_j & if $j\not= i, i+1$\cr
r\tau( {\bf y_{i+1}} \sigma_i^{-1}) f_{i+1}
    + r (1-\tau({\bf y_i}))\tau( \sigma_i^{-1}) f_{i} & if $j=i$\cr
 r \tau(\sigma_i ^{-1})f_{i} & $j=i+1$.}
\end{equation}
We also have the less-interesting  {\it left} action of ${\bf B}_n$
 on this module which is semilinear for the $B_n$ action by conjugation
on ${\cal M},$ and is given by
\begin{equation}
\label{equation:the less interesting left action on H(.., y,..)}
\sigma_i \cdot r f_j = \tau( \sigma_i) r f_j.
\end{equation}
The homology class of an arc $\delta$ based at $y$ is given in terms
of its homotopy class $g\in {\bf F}_n$ by the derivation $e(w)$ defined as
follows:
\begin{equation}
\label{equation:derivation based at y}
e:{\bf F}_n \to \oplus {\cal M}f_i = H_1(P_n,y,{\cal
M}_{{\rm right}}),
 y_i\mapsto f_i,  e(xy)=xe(y)+e(x).
\end{equation}
and we have $[\delta]_y = e(g)$, 
and the derivation $e$ is equivariant in the sense that
$e({\beta}g{\beta}^{-1})={\beta}e(g){\beta}^{-1}.$ \bs
 
Now that we have defined two modules, 
define  our  ${\cal M}$ bilinear pairing  $\langle \   ,\  \rangle$
$$(\oplus {\cal M}f_i) \times (\oplus e_i {\cal M}) \to {\cal M}$$
 by 
\begin{equation}
\label{equation:intersection pairing}
      \langle f_i,e_j\rangle = \left\{\matrix{0, i \ne j\cr      
                  \tau({\bf x_1...x_i})-\tau({\bf x_1...x_{i-1}})   ,
i=j}\right.
\end{equation}
Since the $\tau({\bf x_i})-1$ are non-zero-divisors they will become invertible
in any flat Artinian fraction ring of ${\cal M}$. In other words
this becomes a perfect pairing, i.e. $\tau$ becomes non-degenerate.
This completes the
definition of the pairing. \bs

{\bf Remarks on notation:} We could have avoided introducing the symbols
$e_i, f_i, d$ and $e$, using instead 
$$e_i = [x_i]_x, \ \ f_i = [y_i]_y, \ \ d({\bf w}) = [w]_x,\ \ e({\bf w}) =
 [w]_y$$ where in the final two expressions $w$ is taken to be a
loop based at $x$ or $y$ respectively. 
However we use the symbols $e_i$ and $f_i$  to make them look like basis
elements in two vector spaces, and using the symbols $d$ and $e$ made them look
 like derivations. For this reason we chose what some readers may view as
an excess of notation. \bs

{\bf Proof of (a):} To verify our conventions, let us show that the
pairing defined in this way really is compatible with the braid action.
Thus for ${\beta}\in {\bf B}_n$ we need to check that
$\langle{\beta}v{\beta}^{-1},{\beta}w{\beta}^{-1}\rangle$ = ${\beta}\langle
v,w\rangle{\beta}^{-1}.$ We may assume that the braid homeomorphism ${\beta}$
is
$\sigma_i$ and that our vectors $v$ and $w $ are of the form
$$v=af_i+bf_{i+1}, \ \ \ \ w=e_ir+e_{i+1}s.$$
We have
$$\langle{\beta}v{\beta}^{-1},{\beta}w{\beta}^{-1}\rangle$$
$$=\langle \sigma_i(af_i+bf_{i+1})\sigma_i^{-1},\sigma_i(e_ir+e_{i+1}s)
\sigma_i^{-1}\rangle$$
$$=\langle\sigma_ia[{\bf y_{i+1}}\sigma_i^{-1}f_{i+1}+(1-{\bf
y_i})\sigma_i^{-1}f_i]
+\sigma_ib[\sigma_i^{-1}f_{i+1}],$$
$$[e_{i+1}\sigma_i]r\sigma_i^{-1}+
[e_i\sigma_i{\bf x_i}+e_{i+1}\sigma_i(1-{\bf x_{i+1}})]s\sigma_i^{-1}\rangle.$$
In the formula above the first two expressions in square brackets
come from the right side of (\ref{equation:the interesting right action on
H(.., y,..)}) and the last two  expressions in square brackets come from the
right side of (\ref{equation:the interesting left action on H(.., x,..)}).
\bs

Now using definition  (\ref{equation:action of Bn on Fn}) from
section \ref{subsection:The group B} and also the
definition of the ${\bf y_i}$ we know
$${\bf y_i}\sigma_i^{-1}=\sigma_i^{-1}{\bf y_{i+1}y_iy_{i+1}^{-1}}$$
$${\bf y_{i+1}}\sigma_i^{-1}=\sigma_i^{-1}{\bf y_i}$$
$$\sigma_i{\bf x_i} ={\bf x_{i+1}}\sigma_i$$
$$\sigma_i{\bf x_{i+1}}={\bf x_{i+1}x_ix_{i+1}^{-1}}\sigma_i.$$
Making these substitutions, and using the notation such that
$^\beta a$ denotes $\beta a \beta^{-1},$ we obtain

$$\langle [^{\sigma_i} a(1-{\bf y_iy_{i+1}y_i^{-1}})
+^{\sigma_i} b]f_i +^{\sigma_i} a {\bf y_i}f_{i+1}\ \ , $$
$$ e_{i+1}[(^{\sigma_i}r + (1-{\bf x_{i+1}^{-1}x_ix_{i+1}})^{\sigma_i}s)]
+ e_i {\bf x_{i+1}}^{\sigma_i}s\rangle$$
and using the definition  (\ref{equation:intersection pairing}) of
the intersection pairing this becomes
$$[^{\sigma_i} a(1-{\bf y_iy_{i+1}y_i^{-1}})
+^{\sigma_i} b]({\bf x_1...x_i-x_1...x_{i-1}})[{\bf x_{i+1}}^{\sigma_i}s]$$
$$+[^{\sigma_i} a {\bf y_i}]({\bf x_1...x_{i+1}}-{\bf x_1...x_i})
[(^{\sigma_i}r + (1-{\bf x_{i+1}^{-1}x_ix_{i+1}})^{\sigma_i}s)]
$$
and when this expression is multiplied out, taking
into account ${\bf y_i}={\bf x_1....x_{i-1}x_i^{-1}...x_1^{-1}}$
everything cancels
except
$$^{\sigma_i}a({\bf x_1...x_{i-1}x_{i+1}}-{\bf x_1...x_{i-1}})^{\sigma_i}r +
^{\sigma_i}b({\bf x_1...x_{i+1}}-{\bf x_1...x_ix_{i+1} })^{\sigma_i}s.
$$
Using $\sigma_i{\bf x_1...x_i}\sigma_i^{-1}={\bf x_1...x_{i-1}x_{i+1}}$
this becomes
$$\sigma_i (a({\bf x_1...x_i}-{\bf x_1...x_{i-1}})r+b({\bf x_1...x_{i+1}}-
{\bf x_1...x_i})s)\sigma_i^{-1}$$
$$=\sigma_i\langle
af_i+bf_{i+1},e_ir+e_{i+1}s\rangle\sigma_i^{-1} =
\sigma_i\langle v,w\rangle\sigma_{i} ^{-1}$$
as desired.
This proves statement (a) of the theorem. \bs

{\bf Proof of (b):}  We calculate, using (a). 
$$\langle v\beta,w\rangle = \beta\beta^{-1}\langle v\beta,w\rangle =
\beta\langle\beta^{-1}v\beta,\beta^{-1}\beta w\beta^{-1}\beta\rangle$$
$$ =
\beta [ \beta^{-1}\langle v,\beta w \beta^{-1}\rangle \beta] = \langle v,\beta
w\rangle$$ 

{\bf Proof of (c):}  Left to the reader.

\

{\bf Proof of (d):} We can use the anti-involution $*$ on $  {\cal M}$
 to relate the homology class of an arc in either of our
two homology groups as  follows: Let $e_i^\star\in \oplus_i {\cal M}f_i$
be defined
$$e_i^\star=e({\bf x_i^{-1}})$$
where $e$ is our derivation. Then $e_i^\star$ is a basis and we have
$$\oplus  e_i{\cal M} \to  \oplus {\cal M}f_i$$
$$ \sum e_ir_i \mapsto \sum r_i^\star e_i^\star, \ \ \ \ r_i \in {\cal M}.$$
The inverse map sends $f_i$ to $f_i^\star =d({\bf y_i^{-1}}).$

\

{\bf Proof of (e):} \ 
We are given $v=[w]_x$, where $w$ is a simple loop based at $x$. Every simple
$x$-based loop in $P_n$ is the image of one of the generators,
 ${\bf x_i},$ of
${\bf F}_n$ under a braid homeomorphism. Therefore $v = [\beta(x_i)]_x=\beta
[x_i]_x\beta^{-1}=\beta e_i\beta^{-1}$   for some $\beta\in {\bf B}_n.$ (We
could
even take $i=1.$) Then we have
$$\langle v^*,v\rangle=\langle \beta e_i^* \beta^{-1},\beta e_i \beta^{-1}
\rangle
=\beta \langle e_i^\star,e_i\rangle \beta^{-1}.$$
Now we have $e_i^\star=e({\bf x_i^{-1}})$ where $e$ is the derivation
sending ${\bf y_i}$ to $f_i.$ Writing 
$${\bf x_i^{-1} = {\bf y_1y_2}\dots{\bf y_iy_{i-1}^{-1}}\dots{\bf y_1^{-1}}}$$
we expand using the rule $e({\bf ab})={\bf a} d({\bf b})+e({\bf a}).$ We find
that the term involving  $f_i$ is ${\bf y_1y_2}\dots{\bf y_{i-1}}f_i$ and
therefore
$$\langle e_i^\star,e_i\rangle={\bf y_1...y_{i-1}}\langle f_i,e_i \rangle$$
$$=\tau({\bf y_1}\dots{\bf y_{i-1}})(\tau({\bf
x_1}\dots{\bf x_i})-\tau({\bf x_1}\dots{\bf x_{i-1}}))$$ which simplifies to
$\tau({\bf x_i})-1.$ Then
$$\langle v^\star, v \rangle = \beta^{-1}(\tau({\bf x_i})-1) \beta
=\tau(\beta({\bf x_i}))-1 = \tau({\bf w}) - 1.$$

{\bf Proof of (f):} \ This part of the theorem relates the
intersection pairing with the  the
geometric intersections of arcs on surfaces. Throughout this proof we are
assuming  that
$\tau^+$ is faithful. Therefore we needn't make any distinction between a
braid and its associated matrix.  We need a braid element
$\omega$ in the kernel of the
representation $\tau$. (This is where we use the hypothesis that $\tau$ is
non-faithful for some $r$.) Suppose our element $\omega$ is a braid on
$r$ strands. Say $\delta$ is a  loop based at
$y$ which encircles only $q_i$ and $\gamma$ is a loop based at $x$ which
encircles only $q_j$ and that $\langle [\delta]_y,[ \gamma]_x\rangle=0.$
Recall $y_k$ is the arc based at $y$ whose homotopy class is
${\bf y_k}.$ We may find a braid $\kappa$ such that $\kappa(\delta)=y_k$ for
some number $k$ and we have
$$\langle [\delta]_y,[\gamma]_x\rangle=\kappa\langle
f_k,\kappa^{-1}[\gamma]\kappa \rangle \kappa^{-1}
=\kappa \langle f_k,[\kappa^{-1}\gamma]\rangle \kappa^{-1}.$$
Let us re-name our arcs,
giving $\kappa^{-1}\gamma$ the name $\mu.$ Now we are in the
situation where
$$\langle f_k,[\mu]_x\rangle=0$$
and we wish to show $\mu$ can be homotoped to not meet the arc $y_k.$
First note that $\mu$ does not encircle $q_k$
because if we expand $[\mu]$ in the form $e_1c_1+...+e_nc_n$
then $0=\langle f_k,\mu\rangle=\tau({\bf x_1...x_{i-1}})(\tau({\bf x_i})-1)c_k$
forces $c_k$ to be zero, because $\tau({\bf x_i})-1$ is a nonzero divisor.
But $c_k$ would be a sum of  elements of the form
$\pm\tau({\bf w}) $ for ${\bf w}\in F_n$ whose coefficient
sum adds to the integer
$+1$ if $\mu$ encircled $q_k$ and
such a sum is congruent modulo the $\tau({\bf x_i})-1$ to the element
$1\in {\cal M}$ . By hypothesis the ideal generated by the
$\tau({\bf x_i})-1$ in the appropriate subring
is not the unit ideal so this is impossible.
By enlarging $n$, in effect adding $r-1$ new points numbered
$q_{k+1},...,q_{k+r-1}$ to our
set $Q_n$ and renumbering, 
we may assume $\mu$ does not meet the arcs $y_{k+1},...,y_{k+r-1}$
 based at $y$ whose
homotopy classes are ${\bf y_{k+1},...,y_{k+r-1}},$ so we may assume
$$\langle f_i,[\mu]_x\rangle=0$$ 
for $i=k+1,...,k+r-1.$ By further enlarging $\{q_i\},$ adding
new points $q_1,...,q_r$ and 
renumbering, we may assume
$\mu$ also does not meet the arcs $y_1,...,y_r$ based at $y$ whose homotopy
classes are ${\bf y_1,...,y_r}$ so we also have
$$\langle f_i,[\mu]_x\rangle=0$$
for the values $i=1,2,...,r$ as well.  In all, we
have $\langle f_i, [\mu]_x \rangle=0$ for $i=1,2,3,...,r, k, k+1,...,k+r-1.$  \bs

We  can let $\sigma$ be a second copy of the braid $\omega$ 
but which acts on strands $k,k+1,...,k+r-1$ 
so that the left action of $\sigma$
on only the basis elements $e_k,...,e_{k+r-1}$
is nontrivial, and the action of $\sigma$ on ${\cal M}$ is trivial. \\ \bs

Let $x_r$  as usual be the simple arc based at $x$ 
which encircles just $q_r,$ whose homotopy class is ${\bf x_r}.$
Let $\alpha$ be a braid which does not affect strands
1,2,...,r-1 such that $\alpha(x_r)=\mu.$
We claim that the conjugate braid
$\alpha^{-1}\sigma \alpha$ 
fixes the vectors $e_1,..,e_r \in \oplus e_i {\cal M}.$ It clearly
fixes $e_1,...,e_{r-1}$ since the braid does not involve the corresponding
strands. As for $e_r$  write
$$\alpha e_r = e_1c_1 + ... + e_n c_n.$$
The coefficients $c_i$ which are zero are the same as the ones which
are zero if instead we calculated $\alpha e_r\alpha^{-1}.$ Thus we calculate
$$\langle f_i,\alpha e_r \alpha ^{-1}\rangle = \langle f_i,
\alpha [x_r]\alpha^{-1}\rangle 
=\langle f_i,[\alpha(x_r)]\rangle=\langle f_i,[\mu]\rangle$$
which is zero for
$i=k,...,k+r-1.$ This proves $c_i=0 $ for 
$i=k,k+1,...,k+r-1.$ Since the braid $\sigma$
is in the kernel of the action on ${\cal M}$ it acts only on coefficients
which correspond to the strands it acts on, ie only on
coefficients $c_k, c_k+1,...,c_{k+r-1}$ of $\alpha e_r, $  and yet
these coefficients are all zero. This shows
$\sigma$ fixes $\alpha e_r$ so $\alpha^{-1}\sigma \alpha$ 
fixes $e_r$ as claimed.

\

Now one sees that the commutator $[\alpha^{-1}\sigma \alpha,\omega]$
also fixes the vectors
$e_1,...,e_r.$
This is because $\omega$ preserves
$\oplus_{i=1}^r e_i{\cal M}$ and $\alpha^{-1}\sigma \alpha$ has
no effect on this module.  Now write
$$\beta=\alpha^{-1}\sigma \alpha$$
and let us calculate the action of a commutator
$$[\beta^s, \omega^t]$$ on the whole module, for numbers $s$ and $t.$
We already know $e_1,...,e_r$ are fixed. Thus let $i\ge r+1$
and let us calculate. Write
$$\beta^{-s} e_i=e_ra_r +...+e_na_n$$
and
$$\omega^t e_r=e_1d_1+...+e_rd_r.$$
Then
$$\beta^s\omega^t\beta^{-s}w^{-t}e_i=\beta^s
\omega^t\beta^{-s}e_i=\beta^s\omega^t(e_ra_r+...+e_na_n)$$
$$=\beta^s(e_1d_1a_r+...+e_rd_ra_r +e_{r+1}a_{r+1}+...+e_na_n)$$
$$=\beta^s(e_1d_1a_r+...+e_rd_ra_r-e_ra_r+\beta^{-s}e_i)$$
$$= e_1d_1a_r+...+e_r(d_r-1)a_r +e_i.$$
We see that the action of the commutator adds to the $i'th$ basis
vector only some multiples of basis vectors $e_1,...,e_r.$
In other words, the action of the commutator lies in an abelian group
of elementary matrices with entries in ${\cal M}.$ It follows that for
different values of $s$ and $t$  the commutators all commute.
In other words,
$$1=[[\beta^s,\omega^t],[\beta^{s'},\omega^{t'}]]$$
for all $s,t,s', t'.$

\

Any such relation (for any given values of these numbers)
continues to hold if we substitute $w$ and $b$ for any powers, and
it follows from McCarthy's Tits alternative for braid groups \cite{Mc}
 that $\beta$ and $\omega$ commute.

\

Now, we could have chosen $\omega$ pseudo-Anosov, and so $\beta$ and $\omega$
commuting that means $\beta$ is a braid which does not involve the $r$'th braid
 generator $\sigma_r. $  This means that $\beta$ (which we already know does
not involve the braid generators $\sigma_1,\dots,\sigma_{r-1}$) also does not
involve $\sigma_r$.   In other words,
$\alpha^{-1}\sigma
\alpha$ preserves the homotopy type of
$x_r$. Then $\sigma$ preserves the homotopy type of $\alpha(x_r)=\mu$
and again choosing $\sigma$ pseudo-Anosov on strands $k,...,k+r-1$ we
see this means $\mu$ does not meet the arc with homotopy
class $y_k,$ as needed. This completes the proof of (f) and so also of
Theorem \ref{theorem:intersection pairing}. $\|$

\section{Detecting reducing loops}
\label{section:detecting reducing moves}
We are now ready to make the connection between the
existence of a reducing loop in a closed braid
$\tilde{\beta}$ and our intersection pairing on arcs in the n-times punctured
plane
$P_n = P - Q_n$, under the action of
$\beta\in {\bf B}_n$:

\begin{theorem}
\label{theorem:detecting reducibility}
A braid $\beta \in {\bf B}_n$  is conjugate to a positively (respectively
negatively) reducible braid if and only if there is a simple homology
class $v\in H_1(P_n,x,{\cal M}_{{\rm left}})$
such that $\langle v^\star,{\beta} v {\beta}^{-1}\rangle = 0$ (resp.
 $\langle \beta v^\star \beta^{-1}, v \rangle=0$).
\end{theorem}

\pf  We will treat the positively reducible case, the negative case
being similar. Suppose we have a simple homology class $v$
such that $\langle v^\star,{\beta} v {\beta}^{-1}\rangle=0 $.  The homology
class
$v$, being assumed simple, is represented by an arc $\eta$ based at
$x$ encircling some
$q_i.$   We know that the homology class $v=[\eta]_x$ satisfies
\begin{equation}
0=\langle v^\star,{\beta} v {\beta}^{-1}\rangle= \langle [\eta]^\star,
[{\beta}\eta]\rangle.
\end{equation}
Suppose  that $\eta $ encircles
$q_i $ and ${\beta}\eta$ encircles $q_j.$ We wish to apply  Theorem
\ref{theorem:intersection pairing}. Let $\omega$ be a loop
based at $y$  which represents
the same element of $\pi_1(P_n,y)$ as $\eta$. Then we have
$0=\langle [\omega]_y, [\beta \eta]_x\rangle $ and
the theorem says we can
homotop the two loops $\omega$ and $\beta(\eta)$
apart while fixing both baseoints.
Another way of describing the same information is that if we had
had connected the basepoint $x$ to $q_i$ by an arc $\gamma$ entirely contained
within the closed loop $\eta$ then $\gamma \cup \beta(\gamma)$
could be  
be  homotoped fixing both endpoints $q_i$ and $q_j$ 
so that the union  is a simple arc from $q_i$
to  $q_j$ passing through the basepoint $x.$  Near the basepoint the
orientation of this arc is clockwise (in the direction from $y$ to $x$)
and this is why it is a positive reducing move.  Moreover, rather than modify
${\beta}\gamma$ by a homotopy,  we can choose the braid isotopy so that
$\gamma$ and ${\beta}\gamma$ actually meet only at
the basepoint. 
Finally, it helps to move the basepoint away from the axis $A$ slightly, to a
a new point $x'$ nearby in the interior of $P.$ 
We will be done if we can prove the following lemma,
which actually shows how the reducing move may be carried out. Note that
the lemma does not distinguish between positive and negative reducing moves.

\begin{lemma}
\label{lemma:} Let ${\beta}$ be a braid homeomorphism.
Let $x' \in P$  be a basepoint which is outside the disc where
${\beta}$ acts. Let $\gamma$ be a simple arc in $P$ which joins $x'$ to  $q_i$. 
Then there is a reducing move replacing the  $i^{th}$ strand of the closed braid
$\tilde{\beta}$  by an arc in the plane
$P$ which joins $q_i$ to ${\beta}(q_i)$,  if and only if $\gamma$  can be
homotoped to intersect
${\beta}(\gamma)$ precisely in the point $\{x'\}$ 
\end{lemma}

\pf Suppose $\gamma$ can be homotoped to meet $\beta{\gamma}$ only
at $x'$. Composing the  $\beta$ with a suitable homeomorphism we 
may assume $\gamma$ actually meets $\beta(\gamma)$  only at $x'.$
Recall the map ${\bf H}:P\times [0,1]\to {\mathbb R}^3$ from
section \ref{subsection:braid homeomorphisms and geometric braids}
Consider the image
 under ${\bf H}$   of the square  $\gamma \times [0,1]$.  The  square
 has four edges: 
$$\gamma \times \{0\}, \gamma \times \{1\}, \ \{x\}  \times 
[0,1] \ \ {\rm  and} \ \  \{q_i\} \times [0,1].$$
The map ${\bf H}$  fixes the arc $\gamma = \gamma \times \{0\}$,
sends $\gamma \times 1$  to  ${\beta}(\gamma)$, sends $\{x'\}\times [0,1]$
to a
small homotopically trivial loop about the axis $A$
 and sends $q_i \times [0,1]$ to the $i^{th}$ strand  of  the closed
braid
associated to ${\beta}$ in ${\mathbb R}^3$. 
 The square describes a homotopy by which the $i^{th}$ strand can be
homotoped,
 fixing both ends (so the rest of the braid is not affected) to the
composite of three paths: ${\beta}(\gamma)$,  a small homotopically trivial
closed
loop  about the axis $A$, and $\gamma$.
You can cut out the homotopically trivial closed loop without affecting
the homotopy type of the link in ${\mathbb R}^3$,  and  what
remains is the original closed braid in ${\mathbb R}^3$ with its $i^{th}$
strand
replaced by the arc
$\gamma\cup {\beta}(\gamma)$  in the half  plane $P$. This is the desired reducing move.  The
converse should also be clear, i.e.
if a reducing move replaces the $i^{th}$ strand by 
 an arc in $P$ joining $q_i$ to
${\beta}(q_i) = q_j$, then the arc can always be deformed, in the complement of
$Q_n$,
to an arc of the form $\gamma\cup {\beta}(\gamma)$, where $\gamma$ runs from
$x'$ to $q_i$. This completes the proof of the Lemma, and so also of the
Theorem.
\endpf

{\bf Example:}  We use Theorem \ref{theorem:detecting reducibility} to show
that the braid  $\beta_2 = \sigma_2^{-2}
\sigma_1^{-1}\sigma_2^{-1}\sigma_3^{-1}\sigma_2^3\sigma_1\sigma_2 \sigma_3$
 which is illustrated  in Figure \ref{figure:reducible1}(a) has a
negative reducing loop.  It will be clear later why we have named this
braid $\beta_2.$  We do this by exhibiting a simple homology class
$v\in H_1(P_n,x,{\cal M}_{{\rm left}})$
such that  $\langle \beta_2 v^\star \beta_2^{-1}, v \rangle=0$.  We claim that
we can take $v= [x_3]_x = e_3$.  To see this, use the action given in
(\ref{equation:the interesting left action on H(.., x,..)}) to verify that
$\beta_2 e_3 = e_1\tau(\beta_2) = e_1\beta_2$, so
that
$\langle \beta_2 e_3^\star \beta_2^{-1},e_3\rangle = \langle f_1,
e_3\rangle$, which is indeed zero.  \bs

{\bf Remark 1:} The algebraic action of the braid $\beta_2= \sigma_2^{-2}
\sigma_1^{-1}\sigma_2^{-1}\sigma_3^{-1}\sigma_2^3\sigma_1\sigma_2 \sigma_3$ on
$\pi_1({\bf F_n})$ was given in (\ref{equation:action
of Bn on Fn}). This algebraic action is right to left
(functional notation), so that $\beta_2 {\bf x_3} \beta_2^{-1} = {\bf
x_1}$. This  is opposite to the geometric action which is we depicted earlier,
in Figure
\ref{figure:reducible1}(d). The arc we called $\delta^\prime$
determines the loop $x_1$. We obtain from it the arc we call
$\beta(\delta^\prime)$ (which determines $x_3$) by pushing it one
full turn around the oriented braid axis, keeping its endpoint on $\tilde{K}$.
Thus the algebraic action is opposite in sense to the geometric action.
\bs

{\bf Remark 2:}  It should be clear to the reader, from this example, that the
question of enumerating simple homology classes $e_1a_1 + e_2a_2+\cdots +
e_na_n$ in $H_1(P_n,x,{\cal M}_{{\rm left}})$ is the key to applying Theorem
\ref{theorem:detecting reducibility} to obtain an algorithm for recognizing
reducing loops. See the discussion on open problems at the end of the paper.

\section{Detecting exchange moves}
\label{section:detecting exchange moves} 
Recall the conjugacy class of a closed braid $\tilde{K}$ is said to
{\it admit an exchange move} if  it contains a representative
of the form $P\sigma_{n-1}Q\sigma_{n-1}^{-1}.$

\begin{theorem}
\label{theorem:detecting exchange moves}
 A braid $\beta\in {\bf B}_n$ admits an exchange move
if and only if there are simple homology classes
$v,w \in H_1(P_n,x,{\cal M}_{{\rm left}})$
 such that  $\langle v^\star,w\rangle=0$ and $\langle v^\star,{\beta}
w{\beta}^{-1}\rangle=0. $  
\end{theorem}

Observe that if,  instead $v$ and $w$ are such that
$\langle v^*,w\rangle=\langle\beta v^* \beta^{-1},w\rangle=0$ then
`conjugating' the 
second formula by $\beta$ gives  
$\langle v^*,w\rangle=\langle v^*,\beta^{-1}w\beta\rangle=0$. It follows that
$v$ and $\beta^{-1}w\beta$ are two simple classes satisfying the
hypothesis, so either version implies the other. \bs

\pf   The proof of Theorem \ref{theorem:detecting exchange moves} is very closely
related
to the proof of Theorem  \ref{theorem:detecting reducibility}.  Assume that
$$\langle v^\star,w\rangle=\langle v^\star,\beta w {\beta}^{-1}\rangle = 0.$$
We must obtain a pair of arcs $\gamma,\delta$ in $P_n,$ with $\gamma$
based at $y$ and $\delta$ based at $x$ such that $\gamma$ meets neither
$\delta$ nor ${\beta}(\delta).$   We have
$\gamma$ ending at say $q_i$ and $\delta$ ending at $q_j$. 
We may move the basepoints $y$ and $x$ in $P$ a little away from the axis
but still outside the range where $H$ acts.\bs

We may assume that $H(p,t)=p$ for $0\le t \le 1/2,$ and now we do two
homotopies. By doing a homotopy within the square
${\bf H}(\gamma \times [0,1/2]),$
we replace strand $i$ for $0\le t \le 1/2$  by the image under ${\bf H}$ of
$$(\gamma \times \{0\}) \cup (\{q\} \times [0,1/2])\cup
(\gamma \times \{1/2\}).$$

Similarly, we replace strand $j$
for $1/2 \le t \le 1$ by  the image under ${\bf H}$ of
$$(\tau \times \{1/2\} )\cup ( \{r\} \times [1/2,1])
\cup ({\beta}(\tau)\times
\{1\}).$$
The homotopy in the second case takes place within the square
${\bf H}(\tau \times [1/2,1])$.
Now our braid contains a pair of small semicircles about the axis $A,$
and we can `interchange' the positions
of the semicircles in ${\mathbb R}^3$ without affecting
the type of the closed braid.  The homotopies are then reversed and
the result is an exchange move.  \bs

{\bf Remark 3.}
An explicit way to calculate the new braid after the exchange has
taken place is as follows: Assume that  $\langle v^*,w\rangle=0$ and $\langle
v^*, \beta w \beta^{-1} \rangle=0$. Choose a braid  $\psi$ such that
$v=[\psi x_{n-1}]_x$ and
$w=[\psi x_n]_x.$  Choose a braid $\phi$ such that
$v=[\phi x_{n-1}]_x $ and $ \beta w \beta^{-1}=[\phi x_n]_x.$ Then $\beta$ will
be replaced by $\phi \sigma_{n-1}^{-2}\phi ^{-1}\beta \psi
\sigma_{n-1}^2\psi^{-1}.$
\bs

{\bf Example:} An example is the closed 4-braid
$\beta=(\sigma_2^{-2}\sigma_1\sigma_2^{-1})(\sigma_3)
(\sigma_2^3\sigma_1^{-1}\sigma_2)(\sigma_3^{-1}) = P\sigma_3Q\sigma_3^{-1}$. See Figure
6. It's an interesting
example because the associated closed braid represents the unknot, but the
closed braid does not have a reducing loop. However after two exchange moves
(see the right sketch) it is changed to the braid $K$ which we showed earlier
in
Figure \ref{figure:reducible1}. That braid has a positive reducing loop. The
example was discovered by  Hugh Morton \cite{Mor}.  \bs

We now show that the two exchange moves are  detected by the algebra:
\begin{itemize}
\item \underline{The first exchange move:} We need to find
simple homology classes
$v,w \in H_1(P_n,x,{\cal M}_{{\rm left}})$
 such that  $\langle v^\star,w\rangle=0$ and $\langle v^\star,{\beta}
w{\beta}^{-1}\rangle=0. $  We claim that $v = [x_1]_x = e_1$ and
$$w =   [x_2 x_4 x_2^{-1}]_x=d({\bf x_2x_4x_2^{-1}})
=$$
$$d({\bf x_2 }){\bf x_4x_2^{-1}}+d({\bf x_4}){\bf x_2^{-1}}+d({\bf x_2^{-1}})
= e_2 ({\bf x_4}-1){\bf x_2^{-1}}  + e_4 {\bf x_2^{-1}}$$ 
do the job. 
We have $\langle v^\star,w\rangle = \langle e_1^\star,w\rangle = \langle
f_1,w\rangle = 0$ because the coefficient of
$e_1$ in the expression for $w$ is zero. This is the first part of the test for
the first exchange move. Next, a calculation  gives
$$ \beta w  \beta^{-1}  = e_2 
({\bf x_3x_2x_3^{-1}x_2^{-1}}+{\bf x_3^{-1}x_2^{-1}-x_2^{-1}} )
+e_3 ({\bf x_2x_3^{-1}x_2^{-1}} - {\bf x_3^{-1}x_2^{-1}})$$
and this does not involve $e_1$ so that
$$ \langle v^*, \beta w \beta^{-1}\rangle = \langle e_1^*, \beta w \beta^{-1}
\rangle = \langle f_1, \beta w \beta^{-1}  \rangle = 0. $$
This shows that there is an exchange move relating the homology
classes $v$ and $w.$ If one  
calculates the exchange move by the remark above, one sees
that  the braid $\beta$ has been changed to
$\beta_1 =\sigma_2^{-2}\sigma_1^{-1}\sigma_2^{-1}\sigma_3
\sigma_2^3\sigma_1\sigma_2\sigma_3^{-1}$

\item \underline{The second exchange move:} Start with the braid
$\beta_1$.  
 We consider the homology classes
$ v = [ x_1 x_2 x_3 x_4 x_3^{-1} x_2^{-1} x_1^{-1}]_x$ and
$w = [x_1 x_2 x_3 x_2^{-1} x_1^{-1}]_x$  
Note that  $ v^* = f_4$  and   
the rule for applying the derivation $d$ gives 
$$w=d({\bf x_1x_2x_3x_2^{-1}x_1^{-1}}) =$$
$$e_1 ({\bf x_2x_3x_2^{-1}x_1^{-1} - x_1^{-1}})
+ e_2({\bf x_3x_2^{-1}x_1^{-1} - x_2^{-1} x_1^{-1}})
+e_3 {\bf x_2^{-1} x_1^{-1}}$$
 Since the coefficient of $e_4$ in $w$ is zero, whereas  $ v^* = f_4$, we have
$ \langle v^*, w\rangle = 0.$

Now we calculate
                 $\langle v^*, \beta_1 w \beta_1^{-1}\rangle$
and this can be done most easily by noting that in the action of the
braid automorphism $\beta_1$ we have
   $$\beta_1({\bf x_1x_2x_3 x_2^{-1}x_1^{-1}}) =
{\bf x_1x_2x_3x_1x_3^{-1}x_2^{-1}x_1^{-1}}$$ 
so that $$\beta_1 w \beta_1^{-1} =
[x_1x_2x_3 x_1 x_3^{-1} x_2^{-1} x_1^{-1}]_x.$$
When we calculate this using the derivation rule we see that
no basis element $e_4$ occurs, because the word does not involve $x_4$ at all.
Thus $$ \langle v^*, \beta_1 w \beta_1^{-1} \rangle = 
\langle f_4,
\beta_1 w \beta_1^{-1}\rangle = 0.$$  When we calculate the new braid we
obtain a braid $\beta_2$ which differs from $\beta_1$ in
that the $\sigma_3$ crossings have been reversed. The result $\beta_2$ is the
braid in figure 1, and as we have already seen in the example in section
4, this braid satisfies $\langle \beta_2 e_3^*
\beta_2^{-1},e_3\rangle=0$ so it admits a reducing move to a braid on
three strands only.
\end{itemize}

\begin{figure}[htpb!]
\label{figure:exchange3}
\begin{center}
\centerline{\includegraphics[scale=0.80, bb=86 170 487 652]{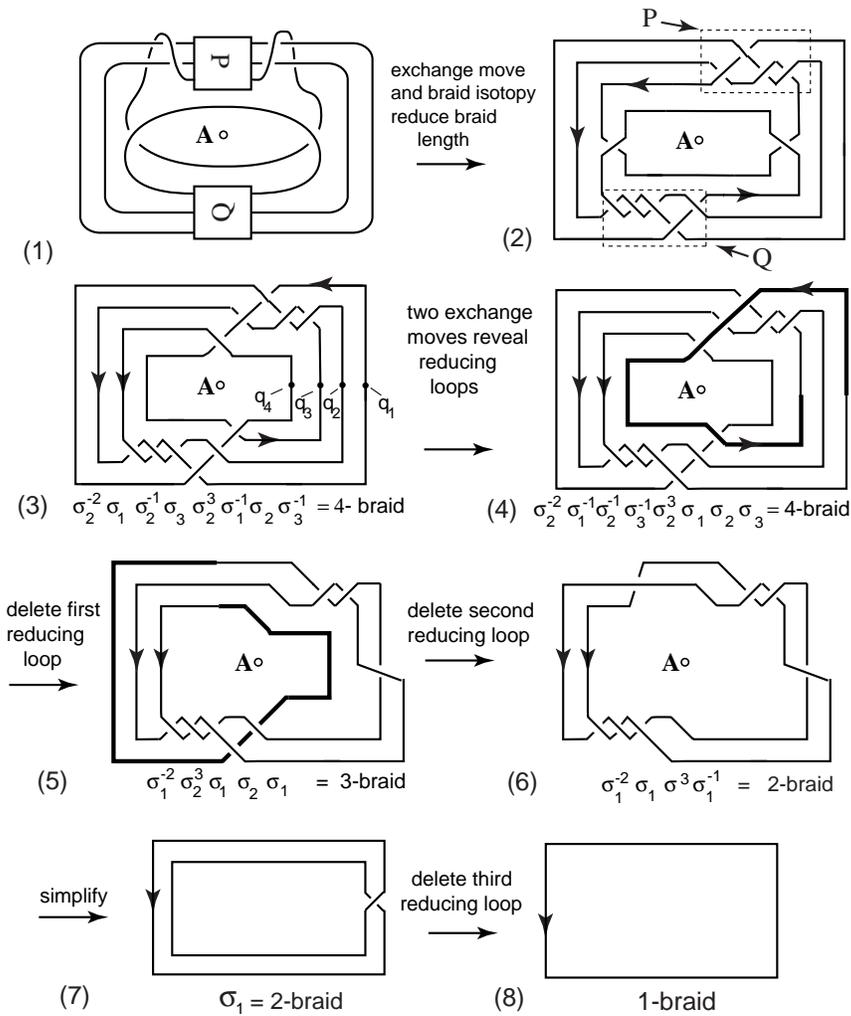}}
\caption{In this closed braid representative of the unknot two exchange
moves are needed before the reducing loop on the right (which is the same as
that in the left sketch in Figure \ref{figure:reducible1}) is revealed.}
\end{center}
\end{figure}

\section{Detecting intersections via  matrix entries}
\label{section:detecting intersections via the monodromy matrices}

In this section we'll explain exactly what it means for a matrix
entry to be zero, and  look at the corresponding very special cases
of the reduction
and  exchange
move conditions. These are weaker than Theorems \ref{theorem:detecting
reducibility} and \ref{theorem:detecting exchange moves}.
Using the notation set up above, with the basis $e_1,...,e_n$   of
$\oplus e_i {\cal M}$ and basis $f_1,...,f_n$   of  $\oplus  {\cal M} f_i$,
recall the intersection pairing which we defined in Equation
(\ref{equation:intersection pairing}). Let $t_i \in {\cal M}$  be
the  element
          $$t_i =    \tau({\bf x_1...x_i}) - \tau ({\bf x_1...x_{i-1}})$$
(Here we use the representation $\tau: {\bf B}_{1,n}\to {\cal M}$ as a method
to
view any element of ${\bf B}_n$ or
${\bf F}_n$ as belonging to ${\cal M}$). \bs

If $\beta\in{\bf B}_n$ is a braid then the image of $\beta$ under our matrix
representation of ${\bf B}_n$ is the matrix  $(r_{ij})\in GL_n({\cal M})$
defined by:
\begin{equation}
          \beta e_j = e_1r_{1j} + ... + e_n r_{nj}.
\end{equation}
First we shall explain  exactly what it means for a matrix
entry $r_{ij}$ to be zero.
Recall if we take {\cal M} to be the Magnus matrix ring then 
each element of ${\cal M}$ is an $n+1$ by $n+1$ matrix with entries in
$\ints[q,q^{-1},t,t^{-1}]$,  so that $(r_{ij})$ can be thought of as
a matrix whose entries are $n+1$ by $n+1$ blocks.\\ \bs 

\begin{lemma}
The matrix entry
 $r_{ij}$ is zero if and only if $\beta=bc$ where strand $i$ at the top
only makes undercrossings in $b$ and strand $j$ at the bottom only makes
overcrossings in $c$, and these are two distinct strands.
 In other words if and only if
$\beta=\sigma_{i+1}\sigma_{i+2}...\sigma_{n-1}P\sigma_{n-1}^{-1} 
Q\sigma_{n-1}\sigma_{n-2}...\sigma_{j+1}$ for some $P,Q\in {\bf B}_{n-1}.$
\end{lemma} 
\pf
We calculate
\begin{equation}
  \langle f_i, \beta e_j\rangle =  t_i r_{ij}
\end{equation}
If $r_{ij}=0$ this means that
\begin{equation}
\label{equation:int2}
           \langle f_i, \beta e_j\rangle=0
\end{equation}
Since $t_i$ is a nonzero divisor 
it follows that  (\ref{equation:int2}) is equivalent
to $r_{ij}=0$.\bs

As usual 
let $x_i$ be the loop at the basepoint $x$  which encircles $q_j$ in
the simplest possible way, so that $e_j = [x_j]_x$ is the homology class of
$x_j$.  From the formulas given earlier we have: 
$$ [{\beta} x_j]_x = {\beta}[x_j]_x{\beta}^{-1}$$
Thus we have
\begin{equation}
\label{equation:int3}
\langle f_i, [{\beta}x_j]_x\rangle
= \langle f_i, {\beta} e_j {\beta}^{-1}\rangle
=    \langle f_i, \beta e_j \tau (\beta)^{-1}\rangle
                    =\langle f_i, \beta e_j\rangle\tau(\beta)^{-1} 
=t_i r_{ij}\tau(\beta)^{-1}
\end{equation}
Since $t_i$ and $\tau(\beta)$
are not zero divisors elements of ${\cal M}$ we see $r_{ij}=0$
if and only if
\begin{equation}
\label{equation:int4}
               \langle f_i, [{\beta} x_j]_x\rangle = 0
\end{equation}
Now $f_i$ itself is the homology class of an arc $y_i$ based at $y$ whose
homotopy class is ${\bf y_i}$. This is an arc which encircles
only $q_i$ but it goes around to the left, and comes in from the
top, and is oriented clockwise.
Also the arc $x_j$ is based at $x$ and encircles $q_j$
anticlockwise. The reason for this difference in orientation is that the basis
elements of $\oplus {\cal M} f_i$ (resp. $\oplus e_i {\cal M}$) are the
$f_i's$ (resp. $e_i's$), which correspond to clockwise (resp anticlockwise)
loops based at $y$ (resp. $x$).
\\
\bs

Now since $f_i = [y_i]_y$ we see from Theorem \ref{theorem:intersection
pairing}
that  (\ref{equation:int4}) is equivalent to the assertion that $y_i$ can be
homotoped  off ${\beta}(x_j).$ \bs

Since ${\beta}(x_j)$ does not now meet $y_i$ then
the word ${\beta}({\bf x_j})$ in the free group does not involve
the letter ${\bf x_i}$. That means
the braid $\beta$ can be factorized bc so that strand $j$ at the bottom only
has overcrossings in $c$ and strand $i$ at the top only undercrossings in
$b$. To see this, consider
 $\beta$ applied to straight arcs from a basepoint $r$ to the $q_s.$
The motions of the $q_s$ are approximated by letting
the $q_s$ move out along these arcs from the basepoint to their final
positions. Now parametrize the $j$'th arc slower than the others, so
in the early part of the braid $q_{j}$ is in front of all other points
(nearer the basepoint). This gives the $c$ part of the braid. And
for a second interval of time only $q_{j}$ moves
and this never moves behind $q_i.$ This gives the $b$ part of the braid.
\endpf\\

Now we can give the corollary describing how to recognize exchange moves
and reducing moves by looking at a single matrix entry. Of course this
is much weaker than Theorems \ref{theorem:detecting reducibility} and
\ref{theorem:detecting exchange moves} because there we find
reducing moves and exchange moves {\it anywhere in the conjugacy class},
but we include this result as it was the first one which we noticed.
\begin{corollary}
\label{corollary:recognizing exchange moves}
\begin{enumerate}
\item
The matrix entry $ r_{n,n-1}$ is zero if and only if
$\beta =P\sigma_{n-1}^{-1}Q\sigma_{n-1}$ for some $P,Q\in {\bf B}_{n-1}.$
In this case $\beta$ admits an exchange move which replaces it
by $Q\sigma_{n-1}^{-1}P\sigma_{n-1}.$
\item
The matrix entry $r_{n,n}$ is zero
if and only if $\beta=P\sigma_{n-1}^{-1}Q$ for some $P,Q\in {\bf B}_{n-1}.$ In
this case $\beta$ admits a reduction move replacing it by $QP.$
\end{enumerate}\end{corollary}
\pf Apply the lemma first in the special case $j=n-1$ and $i=n.$  This
shows $\beta$
is of the form $P\sigma_{n-1}^{-1}Q\sigma_{n-1}.$ 
If $i$ and $j$ are both equal to $n$ then $\beta=P\sigma_{n-1}^{-1}Q.$, and
this
proves the corollary.
 \endpf

But Theorems \ref{theorem:detecting reducibility} and 
\ref{theorem:detecting exchange moves} show there is nothing
special about the loops $\gamma$ and $\delta$ in the proof
of the lemma. Any pair of disjoint simple loops would
give rise to an exchange if the transform of one (based at $x$)
doesn't meet the other (based at $y$). 
Thus, what we have accomplished is to  generalize a
sort of naive notion about matrix entries into a result giving effective
obstructions for reduction or exchange. 

\section{Questions, Comments and Conjectures}
\label{section:questions, comments and conjectures}
We end this paper by discussing the things we have {\it not} been able to do,
mentioning some conjectures and open problems.

\be
\item  A conjecture: Theorem \ref{theorem:intersection pairing}(e) gives a
necessary condition for a homology class $[w]_x \in H_1(P_n,x,{\cal M})$ to be
represented by a simple arc. We conjecture that there is a   series
of finite-dimensional representations $\tau:B_{1,n}\to {\cal M}$
satisfying the conclusions (a) through (f) 
of the Theorem 1, for which this condition is sufficient
as well necessary.  It is even possible that the condition is sufficient in the
special case which we discussed in Section \ref{subsection:Magnus
representations
and an example}. 

\item An important question: Theorem 1 imposes necessary conditions which an
element
$\beta \in Gl_n( {\cal M})$ must satisfy if $\beta$ is to represent a
braid. Namely, for any homology class $v \in \oplus e_i {\cal M}$ we define
$\beta v$ by the ordinary action, and we can define $v^\star \beta$ by the
rule $$v^\star \beta = (\beta^{-1} v)^\star.$$ Then in order to represent
a braid, $\beta $ must satisfy $\langle v^\star \beta, w\rangle=
\langle v^\star, \beta w\rangle$ for all $v,w.$ Secondly, the special
homology class $u=[x_1x_2...x_n]_x$ must satisfy $\beta u = u\tau(\beta).$
It would be interesting to know what other elementary algebraic
conditions characterize 
the set of matrices which represent braids. The more completely this question
can be answered, the closer one comes to finding an algebraic description of
the set of simple homology classes, which in view of theorems 2 and 3
is all that is needed to
immediately recognize whether a braid admits an exchange or reducing move.

\item A comment: Since $\langle v^*,\beta v \beta^{-1}\rangle=\langle v^*,\beta
v\rangle\beta^{-1},$ 
 the condition which is given for recognizing a reducible
braid, in Theorem
\ref{theorem:detecting reducibility},  is equivalent to the apparently simpler
condition $\langle v^*,\beta v\rangle=0.$
However  the more complicated condition has
the advantage that the homology class $\beta v \beta^{-1}$ is simple if and
only if 
$v$ is, and the theorem is an assertion about \underline{simple} homology
classes, 
so we prefer the version we have stated. 

\item A question: One may construct a representation 
$\tau_0:{\bf B}_{1,n}\to {\cal M}$ which describes the action of $B_{1,n}$
on the Laurent polynomial ring $R={\mathbb Z}[T_1^{\pm 1},...,T_n^{\pm 1}].$
The generator $\sigma_i$ acts by permuting $T_i$ and $T_{i+1}$ while the
generator ${\bf x}_i$ acts by multiplication by $T_i$ and ${\cal M}$ is taken
to be the endomorphism ring of $R$ as a module over the subring of symmetric
polynomials.  The matrix entries $r_{ij}$ of the $\tau_0^+$ representation are
elements of $R$ composed with substitutions which permute the variables. Does
the equation $r_{n,n}=0$ characterise braids of the special reducible form
$P\sigma_{n-1}^{-1}Q$ of section \ref{section:detecting reducing moves}?

\item A comment:  When the ${\cal M}$ are matrix rings ${\cal M}=End(V)$
the bilinear form of Theorem  1 can be defined on the actual representation
$V\oplus V \oplus ... \oplus V$ rather than a sum of matrix rings.
Let $\epsilon \in {\cal M}$ such that $V={\cal M}\epsilon$ occurs as a
principal
left ideal. Then the bilinear form on
$\oplus \epsilon^* {\cal M}f_i\ \  \times \ \  \oplus e_i {\cal M}\epsilon$
\ \ \ \ is defined by the formula
$\langle \epsilon^* v, w\epsilon\rangle=\epsilon^*\langle v,w \rangle
\epsilon.$

\item A comment: The bilinear form we have constructed in this
paper is symmetric
in the sense that $\langle v,w \rangle = \langle w^\star, v^\star \rangle.$

\item A comment: In our survey paper
\cite{BLM}, joint with Darren Long, we stated without proof the following
result
generalizing both
\cite{Moo} and
\cite{LP}:
\begin{lemma}
\label{lemma:blm lemma}
 {\rm ( Lemma 2.6 of \cite{BLM})} Suppose the
representations $\tau$ of $B_n$ are not all
faithful. Then for all but finitely many values of $n$
the representation $\tau^+$ of ${\bf B}_n$ is faithful if and only if
the monodromy representation ${\bf F}_n\subset {\bf B}_{1,n}\buildrel
\tau \over \rightarrow {\cal M}$
defines an effective intersection theory for arcs in $P_n.$
\end{lemma}
 D. Long has proven the lemma for
nonabelian monodromy and cap products \cite{Lo}.
In Theorem 1 we establish the Lemma, interpreting the intersection form
as  intersections of two  arcs in $P_n$ based at two different points $x,y$
on the braid axis.     The non-degeneracy of the
intersection pairing is similar to Deligne's non-degenerate
form on cohomology described in $\cite{Lo}$.

\item A comment: When ${\cal M}$ has an anti-involution $\star$ compatible
with $\tau$, as required by the hypotheses of  Theorem 1,
  it is also true in turn
that the ring $Mat_n({\cal M})$ of matrices over ${\cal M}$
has anti-involution compatible with  $\tau^+.$ This is defined such that
a matrix $B$ is sent to the matrix $B'$ such that $(Bf_i^\star)^\star=
f_iB'$ for all $i.$  The process of passing from $\tau$ to
$\tau^+$ may then be iterated, obtaining $\tau^{++}$ and so-on,
and the conclusions of Theorem 1 are retained at each stage.
Although one   now knows that
the question of faithfulness is resolved at the first stage,
our conjecture 1 above allows the possibility of a further `augmentation'
before it may be expected to be true.
\item  Our final remarks deal with relationship between Theorem
\ref{theorem:intersection pairing}, part (f), and Bigelow's `Key Lemma'
(Lemma 3.1 of \cite{Big}).   Our theorem  applies to any representation
 $\tau$ for which $\tau^+$ is faithful,  whereas Bigelow's original
proof involves counting monomial degrees in one particular representation.  Our
theorem in some sense generalizes Bigelow's Key  Lemma to the  cases where
we specialize the parameters $q$ and/or $T$ to values for which
faithfulness is known. For example, by \cite{Kra2} we could specialize $q$ to
any
real number between 0 and 1.  It might be worthwhile to prove this implication
(i.e.
that Krammer's work implies a generalization of Bigelow's Key Lemma)
more precisely. 

Continuing: we interpret the intersection form as intersections of two
loops on the plane $P$ based at distinct points. The same formula generalizes
to any pair of 'noodles.' For the geometric applications, in the
special case ${\cal M}$ is taken to be the Magnus
representation, it would suffice to directly use Bigelow's
lemma in place of our part (f) by replacing $n$ by $n+1$
and approximating one of our arcs by a noodle and the other
by a fork leading to the extra point $q_{n+1}.$ 

Finally, the homology classes we consider are the ordinary
(first) homology classes of our arcs, but with coefficients in
a noncommutative ring. In this theory the self-intersections
of an arc are interesting.
\ee


\begin{thebibliography}{{\bf BM}9}
\markright{} 
\baselineskip18pt


\bibitem{Big1} S. Bigelow, {\it The Burau representation of $B_5$ is not
faithful}, Geometry and Topology {\bf 3} (1999), 397-404.
  
\bibitem{Big} S. Bigelow, {\it Braid groups are linear}, J. AMS 14(2001)
471-486

\bibitem{BiMe-IV} Birman,J.S. and Menasco, W.: {\em Studying Links Via Closed
Braids IV: Closed Braid Representatives of Split and  Composite Links},
Inventiones Math, {\bf 102} Fasc. 1 (1990), 115-139. 

\bibitem {BiMe-V} Birman,J.S. and Menaso, W.: {\em Studying Links Via Closed
Braids 
V: Closed Braid Representatives of the Unlink}, Trans AMS, 
{\bf 329} No. 2 (1992) pp. 585-606. 

\bibitem {BiMe-VI} Birman,J.S. and Menaso, W.: {\em Studying Links Via Closed
Braids VI:A Non-Finiteness Theorem}, Pacific J. Math., {\bf 156}, No. 2,
1992, p. 265-285.

\bibitem{BLMG} J. Birman, {\it Braids, Links and Mapping Class Groups} Annals
of Mathematics Studies 82, Princeton University Press 1974

\bibitem {BH} J. Birman and M. Hirsch, {\it A new algorithm for recognizing the
unknot}, Geometry and Topology {\bf 2} (1998) 175-220

\bibitem {BLMcC} J. Birman, A. Lubotzky and J. McCarthy, {\it Abelian and
solvable
subgroups of the mapping class group}, Duke Math. J. {\bf 50} No. 4 (1983),
1107-1120.

\bibitem {BLM} J. Birman, D. Long, J. Moody, {\it Finite-dimensional
representations
of Artin's braid group}, The Mathematical Legacy of Wilhelm Magnus,
Contemporary Math
{\bf 169} (1994), 123-132.

\bibitem{Bu} W. Burau, {\it Uber Zopfgruppen und gleichsinning verdrillte
Verkettunger}, Abh. Math. Sem. Hanischen Univ. {\bf 11} (1936) 171-178.

\bibitem{FLP} A. Fathi, F. Laudenbach and V. Peonaru, {\it Traveaux de Thurston
sur les surfaces} Asterisque 66-67 (1979), Societe Mathematique de France.

\bibitem{FL} J. Fehrenbach, "Quelques aspects g\'eom\'etriques et dynamiques
du mapping class group" January 1998, PhD thesis, Universit\'e de Nice.

\bibitem{Fiedler} T. Fiedler, {\it A small state sum for knots}, Topology {\bf
32} (1993), 281-294.

\bibitem{Ga} B. J. Gassner, {\it On braid groups}, Abh, Math, Sem. Hamburg
Univ.
{\bf 25} (1961), 19-22.

\bibitem{Iv} N. Ivanov, {\it Subgroups of the Teichmuller Modular Group},
American
Math Soc. Translations of Mathematical Monographs {\bf 115}

\bibitem{Jo} V.F.R. Jones, {\it Hecke algebra representations of braid
groups and link polynomials}, Annals of Math. {\bf 126} (1987), 335-388.

\bibitem {Kra1} D. Krammer, {\it The braid group $B_4$ is linear}, Inventiones
142,
451-48

\bibitem {Kra2} D. Krammer, {\it Braid groups are linear}, Annals of
Mathematics,
to appear.

\bibitem {Law} R. J. Lawrence, {\it Homological representations of the Hecke
algebra}, Comm. Math. Physics {\bf 135} (1990), No. 1,141-191.

\bibitem {Lo} D. Long, {\it Constructing representations of braid groups}, Comm
Math Anal Geom 2(1994) 217-238

\bibitem{LP} D.Long and M. Patton, {\it The Burau representation is not
faithful for
n$\geq$ 6},  Topology {\bf 32} (1993), no. 2, 439-447.

\bibitem {Mc} J. McCarthy, {\it A Tits alternative for subgroups of  mapping
class groups}, Trans AMS 291 (1985) 583-612.

\bibitem {McCool} J. McCool, {\it On reducible braids}, in ``Word Problems II",
editors Adian, Boone and Higman, North-Holland (1980), 261-295.

\bibitem{Men} W. W. Menasco, {\em On iterated torus knots and
transversal knots},{\em Geometry and Topology} {\bf 5} (2001), 651-682.

\bibitem {Moo} J. Moody, {\it The faithfulness question for the Burau
representation}, AMS Proc 119, 1993, 671-679

\bibitem {Mor} H. Morton, {\em An irreducible 4-braid with unknotted
closure}, Math. Proc. Cambridge Phil. Soc., {\bf 93} (1983), 259-261.



\end{thebibliography}
\end{document}